\documentclass[english]{article}
\usepackage[T1]{fontenc}
\usepackage[latin9]{inputenc}
\usepackage{amsmath}
\usepackage{graphicx}
\usepackage{amssymb}
\usepackage{esint}



\usepackage{babel}

\begin{document}

\title{Generalised root identities for zeta functions of curves over finite
fields}

\author{Richard Stone}
\maketitle
\begin{abstract}
We consider generalised root identities for zeta functions of curves
over finite fields, $\zeta_{k}$, and compare with the corresponding
analysis for the Riemann zeta function. We verify numerically that,
as for $\zeta,$ the $\zeta_{k}$ do satisfy the generalised root
identities and we investigate these in detail for the special cases
of $\mu=0,-1\:\&\:-2$. Unlike for $\zeta,$ however, we show that
in the setting of zeta functions of curves over finite fields the
$\mu=-2$ root identity is consistent with the Riemann hypothesis
(RH) proved by Weil. Comparison of this analysis with the corresponding
calculations for $\zeta$ illuminates the fact that, even though both
$\zeta$ and $\zeta_{k}$ have both Euler and Hadamard product representations,
it is the detailed structure of the counting function, $N(T),$ which
drives the Cesaro computations on the root side of these identities
and thereby determines the implications of the root identities for
RH in each setting. 
\end{abstract}

\section{Introduction}

In {[}1{]} and {[}2{]} it was shown that the Riemann zeta function
satisfies the generalised root identities, namely:

\begin{equation}
\frac{-1}{\Gamma(\mu)}\left(\frac{\textrm{d}}{\textrm{d}s}\right)^{\mu}(\ln(\zeta(s)))|_{s=s_{0}}=\textrm{e}^{i\pi\mu}\sum_{\{s_{0}-roots\, r_{i}\: of\:\zeta\}}\frac{M_{i}}{(s_{0}-r_{i})^{\mu}}\label{eq:GenRtId1}\end{equation}
for all $\Re(s_{0})>1$ and arbitrary $\mu\in\mathbb{R}$ (and hence
arbitrary $\mu\in\mathbb{\mathbb{C}}$ by analytic continuation).
In this setting the derivative side is defined by the Euler product
formula, leading to an explicit expression already convergent for
arbitrary $\mu\in\mathbb{\mathbb{C}}$; and in the case of $\mu=1$
the root identity (\ref{eq:GenRtId1}) is equivalent to the Hadamard
product formula for the closely related function $\xi$ (at least
after {}``renormalisation'' and removal of an obstruction - see
{[}1{]} and {[}3, page 35{]}).

Explicit calculation of the values on the root side at $\mu=0,-1\:\&\:-2$
using generalised Cesaro methods was then performed in {[}1{]} and
led, in the case of $\mu=-2$, to a claimed contradiction of the Riemann
hypothesis (RH).

Given the nature of this claim, it is interesting to consider related
settings in which the same techniques may be applied; one such setting
is the class of zeta functions of curves over finite fields, in which
both Euler and Hadamard product formulas hold, but the corresponding
RH is famously known to be true.

The purpose of this paper is to carry out the analogous root identity
computations for this class of zeta functions. By doing so we are
able to confirm that the root identities approach is both applicable
(these zeta functions satisfy the generalised root identities (\ref{eq:GenRtId1}),
just like $\zeta$) and consistent with the truth of the RH for such
zeta functions.

In turn these computations allow us to isolate partially why the RH
is not contradicted by the $\mu=-2$ root identity in the setting
of curves over finite fields, unlike for $\zeta$. Roughly, we find
that it is not merely the existence of both Hadamard and Euler product
formulas that is relevant, but critically also the detailed breakdown
of the counting function $N(T)$ (counting non-trivial roots with
imaginary part in $[0,T)$) into divergent, oscillatory and decaying
asymptotic pieces ($\breve{N}(T)$, $S(T)$ and $\delta(T)$ respectively).

\subsection{Overview}

In section 2 we recall the general form of a zeta function, $\zeta_{k}$,
over a finite field $k$ and use it to deduce both (a) the general
form of the derivative side $d_{\zeta_{k}}(s_{0},\mu)$ of the root
identities for $\zeta_{k}$ and (b) the location of the roots of $\zeta_{k}$
and hence the form of the sum expressing $r_{\zeta_{k}}(s_{0},\mu)$.

Based on this we then show that confirming the root identities for
general $\zeta_{k}$ can be reduced to verifying a particular sub-identity
in which all the generalised roots (i.e. roots or poles here) lie
in an equally spaced fashion on a vertical line $\Re(s)=\sigma_{0}$
with $0\leq\sigma\leq1$ in the complex plane.

Restricting attention to $\Re(s_{0})>1$ as usual, for $\mu>1$ we
then directly verify this sub-identity numerically using the fact
that the resulting sum on the root side is classically convergent.
And since the roots are equally spaced, we then use the Euler-McLaurin
sum formula to extend this readily also to arbitrary $\mu\in\mathbb{R_{<\textrm{1}}}$
by {}``throwing away'' the divergences which then arise, without
any need for additional Cesaro averaging.

It is thus verified that, for arbitrary zeta functions over finite
fields, the generalised root identities

\[
d_{\zeta_{k}}(s_{0},\mu)=r_{\zeta_{k}}(s_{0},\mu)\qquad\forall\Re(s_{0})>1\]
are indeed satisfied for arbitrary $\mu\in\mathbb{R}$ (and hence
arbitrary $\mu\in\mathbb{\mathbb{C}}$ by analytic continuation).

In section 3 we then turn to considering the cases of the $\mu=0,-1\:\&\:-2$
root identities directly, as per {[}1{]} for $\zeta$. It is readily
shown that $d_{\zeta_{k}}(s_{0},\mu)=0$ as a function of $s_{0}$
whenever $\mu\in\mathbb{Z_{\leq\textrm{0}}}$, just as for $\zeta$.
And by again reducing on the root side to the simpler case of a single
vertical line, $\Re(s)=\sigma_{0}$, with equally spaced roots, we
are able to deduce by generalised Cesaro means that we likewise have

\[
r_{\zeta_{k}}(s_{0},0)=r_{\zeta_{k}}(s_{0},-1)=r_{\zeta_{k}}(s_{0},-2)=0\]
as functions of $s_{0}$. Hence we validate explicitly the $\mu=0,-1\:\&\:-2$
root identities and we see, moreover, that unlike for $\zeta$, the
$\mu=-2$ root identity for $\zeta_{k}$ could not lead to a contradiction
of RH in this setting since the Cesaro results for equi-spaced roots
on a vertical line hold equally whether $\sigma_{0}=\frac{1}{2}$
or $\sigma_{0}\neq\frac{1}{2}$.

To confirm this more directly, however, in section 4 we re-perform
the calculation of the root sides $r_{\zeta_{k}}(s_{0},\mu)$ for
$\mu=0,-1\:\&\:-2$ in a way which directly mimics the calculations
in {[}1, section 4.2{]}. We again find $r_{\zeta_{k}}(s_{0},0)=r_{\zeta_{k}}(s_{0},-1)=0$
and $r_{\zeta_{k}}(s_{0},-2)=X_{\epsilon}$, using the notation of
{[}1{]} where the corresponding formula for $\mu=-2$ was $r_{\zeta}(s_{0},-2)=-\frac{1}{2}+X_{\epsilon}$.
The absence of the $-\frac{1}{2}$ term for the $\zeta_{k}$ case
confirms again that, for the setting of zeta functions of curves over
finite fields, the $\mu=-2$ root identity thus does not contradict
the RH for such zeta functions.

We end the paper with some final observations. Firstly we conclude
section 4 by noting why the $\mu=-2$ root identity, which thus requires
$X_{\epsilon}=0$, cannot in fact be used to actually \textit{deduce}
the RH for such $\zeta_{k}$. Then, after summarising in section 5,
we consider fundamental differences in the nature of the counting
functions, $N(T)$, between the cases of $\zeta_{k}$ and $\zeta$,
and discuss how these differences are driving the divergences in the
calculations for $\mu=0,-1\:\&\:-2$ in the two settings (and thereby
the consequences for RH) despite having both Euler and Hadamard product
formulas in both cases.

\section{The Root Identities for Zeta Functions of Curves over Finite Fields}

We adopt the notation of the exposition in {[}3, Chapter 5{]}; $F_{q}$
is the finite field with $q$ elements, where $q=p^{n},\: p$ prime;
and $k$ is a finitely-generated extension of $F_{q}$ of transcendence
degree $1$, so that $k$ is the algebraic extension of $F_{q}(x)$
generated by $y$ satisfying $G(x,y)=0$ for some irreducible polynomial
$G\in F_{q}[x,y]$. Thus $k$ may be viewed as the field of meromorphic
functions on the curve $C$ defined by $G(x,y)=0$.

The zeta function associated to $k$ and $C$ is then defined by the
Euler product formula

\begin{equation}
\zeta_{k}(s)=\prod_{w\in\Sigma(k)}(1-q^{-f(w)s})^{-1}\label{eq:EulerProduct_k}\end{equation}
which converges classically in the half-plane $\Re(s)>1$; here $\Sigma(k)$
is the set of places of $k$ (including $\infty$) and $f:\Sigma(k)\rightarrow\mathbb{Z_{>\textrm{0}}}$
defines the degree of each $w\in\Sigma(k)$. As for $\zeta$, (\ref{eq:EulerProduct_k})
implies that $\zeta_{k}$ has no poles or roots in $\{s:\Re(s)>1\}$.

For purposes of investigating the generalised root identities for
$\zeta_{k}$, however, the starting point is the fact that $\zeta_{k}$,
analytically continued to all of $\mathbb{C}$, has the form

\begin{equation}
\zeta_{k}(s)=\frac{P_{2g}(q^{-s})}{(1-q^{-s})(1-q^{1-s})}\label{eq:zeta_k_general_form}\end{equation}
where $P_{2g}$ is a polynomial with real coefficients of degree $2g$,
$g$ being the genus of the curve $C$. $P_{2g}$ has the property
that it may be factorised as

\begin{equation}
P_{2g}(u)=\prod_{\lambda\in A}(1-\lambda u)\label{eq:P_2g_factorisation}\end{equation}
where the set $A$ is closed under the mapping $\lambda\rightarrow\frac{q}{\lambda}$.%
\footnote{Here is the only place where we amend the notation of {[}3{]} to use
$\lambda$ rather than $\alpha$ since we wish to reserve $\alpha$
for its usual role in later Cesaro calculations.%
} It follows from this and the functional equation for $\zeta_{k}$

\begin{equation}
\zeta_{k}(1-s)=q^{(1-g)(1-2s)}\zeta_{k}(s)\label{eq:zeta_k_functionalEqn}\end{equation}
that the generalised roots (roots and poles) of $\zeta_{k}$ all lie
in the critical strip $0\leq\Re(s)\leq1$ and, like for $\zeta$,
are mirror-symmetric in both the real axis and the critical line $\Re(s)=\frac{1}{2}$.

The restriction of the generalised roots to the critical strip reflects
that $1\leq|\lambda|\leq q$ for all the factors in $P_{2g}(u)$ in
(\ref{eq:P_2g_factorisation}) and if we add to $A$ the two terms
$\lambda=1$ and $\lambda=q$ corresponding to the two factors in
the denominator in (\ref{eq:zeta_k_general_form}), and thus define

\begin{equation}
\tilde{A}=A\,\cup\,\{1,q\}\label{eq:A_tilde_defn}\end{equation}
then we see that $\zeta_{k}$ may be written in (\ref{eq:zeta_k_general_form})
as

\begin{equation}
\zeta_{k}(s)=\prod_{\lambda\in\tilde{A}}(1-\lambda q^{-s})^{\nu_{\lambda}}\label{eq:zeta_k_general_form2}\end{equation}
where $\nu_{\lambda}=-1$ for $\lambda\in\{1,q\}$ and $\nu_{\lambda}=1$
otherwise.

The RH for zeta functions of curves over finite fields is the claim
that, for any such $\zeta_{k}$,

\begin{equation}
|\lambda|=q^{\frac{1}{2}}\qquad\forall\lambda\in A\label{eq:RH_finite_fields}\end{equation}
or equivalently that all the actual roots (not poles) of $\zeta_{k}$
lie strictly on the critical line $\Re(s)=\frac{1}{2}$.

The expression (\ref{eq:zeta_k_general_form2}) is immediately well-adapted
for analysing the generalised root identities (\ref{eq:GenRtId1})
for $\zeta_{k}$. In order to do this, from this point on we restrict
consideration to $\Re(s_{0})>1$.

On the derivative side we immediately have

\begin{eqnarray}
d_{\zeta_{k}}(s_{0},\mu) & = & \frac{-1}{\Gamma(\mu)}\sum_{\lambda\in\tilde{A}}\nu_{\lambda}\left(\frac{\textrm{d}}{\textrm{d}s}\right)^{\mu}\left(\ln(1-\lambda q^{-s})\right)|_{s=s_{0}}\nonumber \\
\nonumber \\ & = & \frac{1}{\Gamma(\mu)}\sum_{\lambda\in\tilde{A}}\nu_{\lambda}\left(\frac{\textrm{d}}{\textrm{d}s}\right)^{\mu}\left\{ \sum_{n=1}^{\infty}\frac{1}{n}\lambda^{n}q^{-ns}\right\} |_{s=s_{0}}\nonumber \\
\nonumber \\ & = & \frac{1}{\Gamma(\mu)}\sum_{\lambda\in\tilde{A}}\nu_{\lambda}\sum_{n=1}^{\infty}\frac{\lambda^{n}}{n}\,\textrm{e}^{i\pi\mu}n^{\mu}(\ln q)^{\mu}q^{-ns_{0}}\nonumber \\
\nonumber \\ & = & \frac{\textrm{e}^{i\pi\mu}}{\Gamma(\mu)}(\ln q)^{\mu}\sum_{\lambda\in\tilde{A}}\nu_{\lambda}\sum_{n=1}^{\infty}\frac{\lambda^{n}}{n^{1-\mu}}\, q^{-ns_{0}}\label{eq:RootId_zeta_k_DerivSide}\end{eqnarray}
which is convergent for arbitrary $\mu$ since $|\lambda q^{-s_{0}}|<1$.

On the root side, for any given $\lambda\in\tilde{A}$, in (\ref{eq:zeta_k_general_form2})
we clearly have infinitely many generalised roots of $\zeta_{k}$
associated to the factor $(1-\lambda q^{-s})^{\nu_{\lambda}}$; if
we write

\begin{equation}
\lambda=q^{\sigma_{0}}\textrm{e}^{i\theta_{0}}\qquad\textrm{where}\qquad0\leq\sigma_{0}\leq1,\:0\leq\theta_{0}<2\pi\label{eq:lambda_defn}\end{equation}
then these are all on the vertical line $\Re(s_{0})=\sigma_{0}$.
There is a unique such root having imaginary part in $[0,\frac{2\pi}{\ln q})$,
which we shall call the {}``base root'' of $\zeta_{k}$ associated
to this $\lambda$-factor and which we shall denote $r_{0}^{(\lambda)}$
so that

\begin{equation}
r_{0}^{(\lambda)}=\sigma_{0}+i\tau_{0}\quad\textrm{where}\quad\tau_{0}=\frac{\theta_{0}}{\ln q}\label{eq:r0lambda_defn}\end{equation}
The other roots associated to this $\lambda$-factor are equally spaced
up and down this vertical line in the complex plane at intervals of

\begin{equation}
C:=\frac{2\pi}{\ln q}\label{eq:C_defn}\end{equation}
i.e the roots of $\zeta_{k}$ arising from the factor $(1-\lambda q^{-s})^{\nu_{\lambda}}$
in (\ref{eq:zeta_k_general_form2}) are the set

\begin{equation}
R_{\lambda}:=\{r_{j}\}_{j\in\mathbb{Z}}\quad\textrm{where}\quad r_{j}=r_{0}^{(\lambda)}+i\cdot Cj=\sigma_{0}+i(\tau_{0}+Cj)\label{eq:R_lambda_defn}\end{equation}

If $\nu_{\lambda}=-1$ these are the poles (all with multiplicity
$M_{i}=-1)$ associated to the two factors on the denominator in (\ref{eq:zeta_k_general_form})
and lying on the boundary-lines of the critical strip $\Re(s_{0})=0$
and $\Re(s_{0})=1$; if $\nu_{\lambda}=1$ they are the roots (all
with multiplicity $M_{i}=1)$ associated to the factors of $P_{2g}$.

On the root side of the generalised root identities (\ref{eq:GenRtId1})
we then have, on noting $\nu_{\lambda}=M_{i}$ in all cases, that

\begin{eqnarray}
r_{\zeta_{k}}(s_{0},\mu) & = & \textrm{e}^{i\pi\mu}\sum_{\lambda\in\tilde{A}}\nu_{\lambda}\sum_{r_{j}\in R_{\lambda}}\frac{1}{(s_{0}-r_{j})^{\mu}}\nonumber \\
\nonumber \\ & = & \textrm{e}^{i\pi\mu}\sum_{\lambda\in\tilde{A}}\nu_{\lambda}\sum_{j=-\infty}^{\infty}\frac{1}{((s_{0}-r_{0}^{(\lambda)})-i\cdot Cj)^{\mu}}\label{eq:RootId_zeta_k_RootSide}\end{eqnarray}
which is classically convergent for $\Re(\mu)>1$ and can then be
analytically continued to $\Re(\mu)\leq1$.

Confirming whether $\zeta_{k}$ satisfies the generalised root identities
thus consists of testing whether the expressions for $d_{\zeta_{k}}(s_{0},\mu)$
and $r_{\zeta_{k}}(s_{0},\mu)$ in (\ref{eq:RootId_zeta_k_DerivSide})
and (\ref{eq:RootId_zeta_k_RootSide}) agree for arbitrary $\Re(s_{0})>1$
and arbitrary $\mu\in\mathbb{R}$ (and hence arbitrary $\mu\in\mathbb{\mathbb{C}}$
by analytic continuation). It is clear, in turn, from these expressions
that this will be true if the contributions on each side from each
$\lambda$-factor $(1-\lambda q^{-s})^{\nu_{\lambda}}$ agree; i.e.
if

\begin{equation}
d_{\zeta_{k}}^{(\lambda)}(s_{0},\mu)=r_{\zeta_{k}}^{(\lambda)}(s_{0},\mu)\label{eq:RootSubId}\end{equation}
where

\begin{equation}
d_{\zeta_{k}}^{(\lambda)}(s_{0},\mu)=\frac{\textrm{e}^{i\pi\mu}}{\Gamma(\mu)}(\ln q)^{\mu}\nu_{\lambda}\sum_{n=1}^{\infty}\frac{\lambda^{n}}{n^{1-\mu}}\, q^{-ns_{0}}\label{eq:RootSubId_DerivSide}\end{equation}
and

\begin{equation}
r_{\zeta_{k}}^{(\lambda)}(s_{0},\mu)=\textrm{e}^{i\pi\mu}\nu_{\lambda}\sum_{j=-\infty}^{\infty}((s_{0}-r_{0}^{(\lambda)})-i\cdot Cj)^{-\mu}\label{eq:RootSubId_RootSide}\end{equation}
and where $\lambda$ and $r_{0}^{(\lambda)}$ are related by (\ref{eq:lambda_defn})
and (\ref{eq:r0lambda_defn}).

For $\mu\in\mathbb{R_{>\textrm{1}}}$, since both sides are convergent,
this may be checked directly numerically. For $\mu\in\mathbb{R_{\leq\textrm{1}}}$
the Euler-McLaurin sum formula (e.g. as stated in {[}4{]}) may be
used to calculate the analytic continuation on the root side explicitly,
since the roots contributing to $r_{\zeta_{k}}^{(\lambda)}(s_{0},\mu)$
are evenly spaced. We have

\begin{eqnarray}
\sum_{j=1}^{k}((s_{0}-r_{0}^{(\lambda)})-i\cdot Cj)^{-\mu} & = & \left\{ \begin{array}{cc}
\frac{i}{(1-\mu)C}((s_{0}-r_{0}^{(\lambda)})-i\cdot Ck)^{1-\mu}\\
\\+f_{+}(s_{0},\mu)\\
\\+\frac{1}{2}((s_{0}-r_{0}^{(\lambda)})-i\cdot Ck)^{-\mu}\\
\\+\frac{i\mu C}{12}((s_{0}-r_{0}^{(\lambda)})-i\cdot Ck)^{-\mu-1}\\
\\+\frac{i\mu(\mu+1)(\mu+2)C^{3}}{720}((s_{0}-r_{0}^{(\lambda)})-i\cdot Ck)^{-\mu-3}\\
\\+\ldots\end{array}\right\} \nonumber \\
\label{eq:EMcL_RootSideA}\end{eqnarray}
and similarly

\begin{eqnarray}
\sum_{j=1}^{\tilde{k}}((s_{0}-r_{0}^{(\lambda)})+i\cdot Cj)^{-\mu} & = & \left\{ \begin{array}{cc}
-\frac{i}{(1-\mu)C}((s_{0}-r_{0}^{(\lambda)})+i\cdot C\tilde{k})^{1-\mu}\\
\\+f_{-}(s_{0},\mu)\\
\\+\frac{1}{2}((s_{0}-r_{0}^{(\lambda)})+i\cdot C\tilde{k})^{-\mu}\\
\\-\frac{i\mu C}{12}((s_{0}-r_{0}^{(\lambda)})+i\cdot C\tilde{k})^{-\mu-1}\\
\\-\frac{i\mu(\mu+1)(\mu+2)C^{3}}{720}((s_{0}-r_{0}^{(\lambda)})+i\cdot C\tilde{k})^{-\mu-3}\\
\\+\ldots\end{array}\right\} \nonumber \\
\label{eq:EMcL_RootSideB}\end{eqnarray}
and since $f_{\pm}(s_{0},\mu)$ clearly represent the analytic continuations
of these sums from $\mu\in\mathbb{R_{>\textrm{1}}}$ to $\mu\in\mathbb{R_{\leq\textrm{1}}}$,
so in general

\begin{eqnarray}
 &  & r_{\zeta_{k}}^{(\lambda)}(s_{0},\mu)=\underset{k,\tilde{k}\rightarrow\infty}{lim}\textrm{e}^{i\pi\mu}\nu_{\lambda}\nonumber \\
\nonumber \\ & \times & \left\{ \begin{array}{cc}
\sum_{j=-\tilde{k}}^{k}((s_{0}-r_{0}^{(\lambda)})-i\cdot Cj)^{-\mu}-\\
\\(s_{0}-r_{0}^{(\lambda)})^{-\mu}-\\
\\\frac{i}{(1-\mu)C}\left\{ \begin{array}{cc}
((s_{0}-r_{0}^{(\lambda)})-i\cdot Ck)^{1-\mu}-\\
\\((s_{0}-r_{0}^{(\lambda)})+i\cdot C\tilde{k})^{1-\mu}\end{array}\right\} -\\
\\\frac{1}{2}\left\{ ((s_{0}-r_{0}^{(\lambda)})-i\cdot Ck)^{-\mu}+((s_{0}-r_{0}^{(\lambda)})+i\cdot C\tilde{k})^{-\mu}\right\} -\\
\\\frac{i\mu C}{12}\left\{ \begin{array}{cc}
((s_{0}-r_{0}^{(\lambda)})-i\cdot Ck)^{-\mu-1}-\\
\\((s_{0}-r_{0}^{(\lambda)})+i\cdot C\tilde{k})^{-\mu-1}\end{array}\right\} -\\
\\\frac{i\mu(\mu+1)(\mu+2)C^{3}}{720}\left\{ \begin{array}{cc}
((s_{0}-r_{0}^{(\lambda)})-i\cdot Ck)^{-\mu-3}-\\
\\((s_{0}-r_{0}^{(\lambda)})+i\cdot C\tilde{k})^{-\mu-3}\end{array}\right\} \\
\\-\ldots\end{array}\right\} \nonumber \\
\label{eq:EMcL_RootSide_Final}\end{eqnarray}
where the number of divergent terms to be subtracted is dictated by
the value of $\mu$ and we may as well take $\tilde{k}=k$.

The formulae for $d_{\zeta_{k}}^{(\lambda)}(s_{0},\mu)$ and $r_{\zeta_{k}}^{(\lambda)}(s_{0},\mu)$
in (\ref{eq:RootSubId_DerivSide}) and (\ref{eq:EMcL_RootSide_Final})
have been implemented in R-code and their equality checked for a variety
of values of $\lambda,\: s_{0}$ and $\mu\in\mathbb{R}$. For example
the graphs below show the equality of $d_{\zeta_{k}}^{(\lambda)}(s_{0},\mu)$
and $r_{\zeta_{k}}^{(\lambda)}(s_{0},\mu)$ for $s_{0}=5.1238$ and
a variety of $\mu$-values ranging from $\mu=-1.5$ to $\mu=2.6$
in the case when $q=5^{2}$ and $\lambda$ is given by $\sigma_{0}=0.6,\:\tau_{0}=\frac{3\pi}{4}$;
here we have used $k=\tilde{k}=1,000$ roots on the root side in (\ref{eq:EMcL_RootSide_Final})
in calculating the approximation to $r_{\zeta_{k}}^{(\lambda)}(s_{0},\mu)$,
while we have approximated $d_{\zeta_{k}}^{(\lambda)}(s_{0},\mu)$
in (\ref{eq:RootSubId_DerivSide}) using $20$ terms in the sum for
each $\mu$. The R-code generating these graphs is given in Appendix
1.

\includegraphics[scale=0.5]{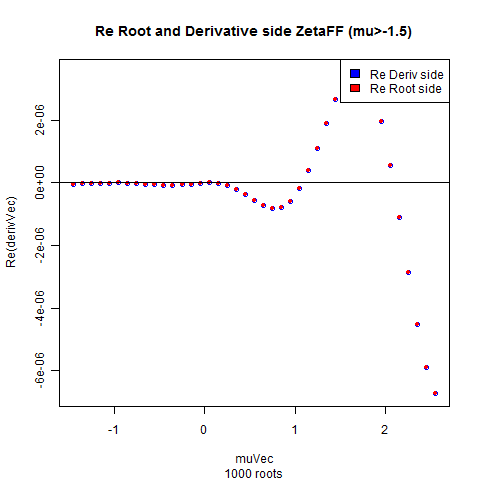}

\includegraphics[scale=0.5]{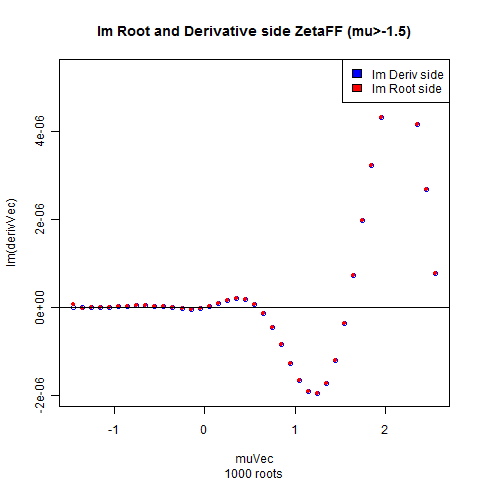}

Based on these sort of numerical verifications we are able to conclude
that:

\paragraph{Result 1:}

\textit{We have}

\begin{equation}
d_{\zeta_{k}}^{(\lambda)}(s_{0},\mu)=r_{\zeta_{k}}^{(\lambda)}(s_{0},\mu)\label{eq:RootId_lambda_factor_true}\end{equation}
\textit{in general for arbitrary} $s_{0}$ \textit{and} $\mu\in\mathbb{R}$
\textit{(and hence arbitrary} $\mu\in\mathbb{C}$ \textit{by analytic
continuation) for any allowed} $\lambda$\textit{-factor. Hence for
any zeta function,} $\zeta_{k}$\textit{, of a curve over a finite
field} $k$ \textit{the generalised root identities are also satisfied,
i.e.}

\begin{equation}
d_{\zeta_{k}}(s_{0},\mu)=r_{\zeta_{k}}(s_{0},\mu)\label{eq:RootId_zeta_k_true}\end{equation}
\textit{for arbitrary} $s_{0}$ \textit{and} $\mu\in\mathbb{R}$ \textit{(and
hence arbitrary} $\mu\in\mathbb{C}$ \textit{by analytic continuation)}.

\section{The root identities for $\zeta_{k}$ at $\mu=0,-1$ and $-2$}

We now turn, as for $\zeta$ in {[}1{]}, to considering the particular
cases of these generalised root identities for $\zeta_{k}$ at $\mu=0,-1$
and $-2$.

On the derivative side in (\ref{eq:RootId_zeta_k_DerivSide}), since
the sum is convergent for arbitrary $\mu$ and the factor $\Gamma(\mu)$
on the denominator diverges for $\mu\in\mathbb{Z}_{\leq0}$, it follows
at once that, as for $\zeta$, we have

\paragraph{Result 2:}

\textit{As a function of} $s_{0}$ \textit{on} $\Re(s_{0})>1$

\begin{equation}
d_{\zeta_{k}}(s_{0},\mu)=0\label{eq:zeta_k_DerivSide_muNegInt}\end{equation}
\textit{whenever} $\mu\in\mathbb{Z_{\leq\textrm{0}}}$.

On the root side, as in section 2, we use the fact that

\begin{equation}
r_{\zeta_{k}}(s_{0},\mu)=\sum_{\lambda\in\tilde{A}}r_{\zeta_{k}}^{(\lambda)}(s_{0},\mu)\label{eq:zeta_k_RootSide_Sum_lambdas}\end{equation}
to reduce the computations to the case of $r_{\zeta_{k}}^{(\lambda)}(s_{0},\mu)$
for any single given $\lambda$-factor $(1-\lambda q^{-s})^{\nu_{\lambda}}$.

In this case the roots are all equi-spaced vertically at intervals
of $C=\frac{2\pi}{\ln q}$ either side of the base root $r_{0}^{(\lambda)}=\sigma_{0}+i\tau_{0}$
on the line $\Re(s)=\sigma_{0}$.

Adopting our standard notation for Cesaro calculations we thus have

\begin{equation}
T=Ck+\tau_{0}+\alpha\qquad\textrm{and}\qquad\tilde{T}=C\tilde{k}-\tau_{0}+\tilde{\alpha}\label{eq:T_Ttilde_defn}\end{equation}
where now $\alpha,\tilde{\alpha}\in[0,C)$; and the variables $z$
and $\tilde{z}$ arising from consideration of $\{s_{0}-roots\: r_{i}\in R_{\lambda}\}$
are given by

\begin{equation}
z=(s_{0}-\sigma_{0})-iT\qquad\textrm{so}\qquad T=i(z-(s_{0}-\sigma_{0}))\label{eq:z_defn}\end{equation}
and

\begin{equation}
\tilde{z}=(s_{0}-\sigma_{0})+i\tilde{T}\qquad\textrm{so}\qquad\tilde{T}=-i(\tilde{z}-(s_{0}-\sigma_{0}))\label{eq:ztilde_defn}\end{equation}

Recalling our definition of Cesaro convergence along a contour in
{[}1, section 2{]} (namely that we (a) first remove generalised Cesaro
eigenfunctions in the \textit{geometric} variables ($z$ or $\tilde{z}$)
and then (b) apply a suitable power of the Cesaro averaging operator,
$P$, to the residual function to reduce to one with a classical limit)
we are then able to deduce the following calculational results which
arise in the computation of $r_{\zeta_{k}}^{(\lambda)}(s_{0},0)$,
$r_{\zeta_{k}}^{(\lambda)}(s_{0},-1)$ and $r_{\zeta_{k}}^{(\lambda)}(s_{0},-2)$:

\paragraph{Lemma:}

(1) \textit{On the contour traversed by} $z=(s_{0}-\sigma_{0})-iT$
\textit{we have}

(a) \[
\underset{z\rightarrow\infty}{Clim}\,\alpha^{n}=\frac{1}{n+1}\, C^{n}\]

(b) \begin{eqnarray*}
\underset{z\rightarrow\infty}{Clim}\, k\alpha & = & -\frac{i}{2}(s_{0}-r_{0}^{(\lambda)})-\frac{1}{3}C\qquad\textrm{and}\\
\\\underset{z\rightarrow\infty}{Clim}\, k\alpha^{2} & = & -\frac{i}{3}C(s_{0}-r_{0}^{(\lambda)})-\frac{1}{4}C^{2}\qquad\textrm{and}\\
\\\underset{z\rightarrow\infty}{Clim}\, k^{2}\alpha & = & -\frac{1}{2C}(s_{0}-r_{0}^{(\lambda)})^{2}+\frac{7}{12}i(s_{0}-r_{0}^{(\lambda)})+\frac{1}{4}C+\frac{1}{12}\tau_{0}\end{eqnarray*}

(c) \begin{eqnarray*}
\underset{z\rightarrow\infty}{Clim}\, z\alpha & = & \underset{z\rightarrow\infty}{Clim}\, z\alpha^{2}=0\qquad\textrm{and}\\
\\\underset{z\rightarrow\infty}{Clim}\, z^{2}\alpha & = & \frac{i}{12}C^{2}(s_{0}-\sigma_{0})\end{eqnarray*}

(d) \begin{eqnarray*}
\underset{z\rightarrow\infty}{Clim}\, k & = & \frac{1}{C}\left\{ -i(s_{0}-r_{0}^{(\lambda)})-\frac{1}{2}C\right\} \qquad\textrm{and}\\
\\\underset{z\rightarrow\infty}{Clim}\, k^{2} & = & \frac{1}{C^{2}}\left\{ -(s_{0}-r_{0}^{(\lambda)})^{2}+iC(s_{0}-r_{0}^{(\lambda)})+\frac{1}{3}C^{2}\right\} \qquad\textrm{and}\\
\\\underset{z\rightarrow\infty}{Clim}\, k^{3} & = & \frac{1}{C^{3}}\left\{ \begin{array}{cc}
\frac{i}{4}C^{2}(s_{0}-\sigma_{0})+i(s_{0}-r_{0}^{(\lambda)})^{3}+\frac{3C}{2}(s_{0}-r_{0}^{(\lambda)})^{2}\\
\\-iC^{2}(s_{0}-r_{0}^{(\lambda)})-\frac{1}{4}C^{3}\end{array}\right\} \end{eqnarray*}
(2) \textit{Similarly, on the contour traversed by} $\tilde{z}=(s_{0}-\sigma_{0})+i\tilde{T}$
\textit{we have}

(a) \[
\underset{\tilde{z}\rightarrow\infty}{Clim}\,\tilde{\alpha}^{n}=\frac{1}{n+1}\, C^{n}\]

(b) \begin{eqnarray*}
\underset{\tilde{z}\rightarrow\infty}{Clim}\,\tilde{k}\tilde{\alpha} & = & \frac{i}{2}(s_{0}-r_{0}^{(\lambda)})-\frac{1}{3}C\qquad\textrm{and}\\
\\\underset{\tilde{z}\rightarrow\infty}{Clim}\,\tilde{k}\tilde{\alpha}^{2} & = & \frac{i}{3}C(s_{0}-r_{0}^{(\lambda)})-\frac{1}{4}C^{2}\qquad\textrm{and}\\
\\\underset{\tilde{z}\rightarrow\infty}{Clim}\,\tilde{k}^{2}\tilde{\alpha} & = & -\frac{1}{2C}(s_{0}-r_{0}^{(\lambda)})^{2}-\frac{7}{12}i(s_{0}-r_{0}^{(\lambda)})+\frac{1}{4}C-\frac{1}{12}\tau_{0}\end{eqnarray*}

(c) \begin{eqnarray*}
\underset{\tilde{z}\rightarrow\infty}{Clim}\,\tilde{z}\tilde{\alpha} & = & \underset{\tilde{z}\rightarrow\infty}{Clim}\,\tilde{z}\tilde{\alpha}^{2}=0\qquad\textrm{and}\\
\\\underset{\tilde{z}\rightarrow\infty}{Clim}\,\tilde{z}^{2}\tilde{\alpha} & = & -\frac{i}{12}C^{2}(s_{0}-\sigma_{0})\end{eqnarray*}

(d) \begin{eqnarray*}
\underset{\tilde{z}\rightarrow\infty}{Clim}\,\tilde{k} & = & \frac{1}{C}\left\{ i(s_{0}-r_{0}^{(\lambda)})-\frac{1}{2}C\right\} \qquad\mathit{\textrm{and}}\\
\\\underset{\tilde{z}\rightarrow\infty}{Clim}\,\tilde{k}^{2} & = & \frac{1}{C^{2}}\left\{ -(s_{0}-r_{0}^{(\lambda)})^{2}-iC(s_{0}-r_{0}^{(\lambda)})+\frac{1}{3}C^{2}\right\} \qquad\mathit{\textrm{and}}\\
\\\underset{\tilde{z}\rightarrow\infty}{Clim}\,\tilde{k}^{3} & = & \frac{1}{C^{3}}\left\{ \begin{array}{cc}
-\frac{i}{4}C^{2}(s_{0}-\sigma_{0})-i(s_{0}-r_{0}^{(\lambda)})^{3}+\frac{3C}{2}(s_{0}-r_{0}^{(\lambda)})^{2}\\
\\+iC^{2}(s_{0}-r_{0}^{(\lambda)})-\frac{1}{4}C^{3}\end{array}\right\} \end{eqnarray*}

\paragraph{Proof:}

We shall prove only the results in part (1) for the contour traversed
by $z$; the adaptation of these arguments for the contour traversed
by $\tilde{z}$ in part (2) is straightforward. 

The result for $\underset{z\rightarrow\infty}{Clim}\,\alpha^{n}$
in (1)(a) is immediate since the function $f(T)=\alpha^{n}$ is periodic
with period $C$ and its average value on each period is $\frac{1}{C}\int_{0}^{C}\alpha^{n}\,\textrm{d}\alpha=\frac{1}{n+1}\, C^{n}$.

For the remaining results in (1)(b)-(1)(d) we then calculate successively
in the following order:

\[
k\rightarrow k\alpha^{n}\rightarrow z\alpha^{n}\rightarrow k^{2}\rightarrow k^{2}\alpha^{n}\rightarrow z^{2}\alpha^{n}\rightarrow k^{3}\]
with the calculation of the Cesaro limit of each relying on the Cesaro
limits of the preceding ones in the chain.

Starting with $k$ we have\begin{eqnarray}
k & = & \frac{1}{C}(T-\tau_{0}-\alpha)=\frac{1}{C}\left\{ i(z-(s_{0}-\sigma_{0}))-\tau_{0}-\alpha\right\} \nonumber \\
\nonumber \\ & = & \frac{1}{C}\left\{ iz-i(s_{0}-r_{0}^{(\lambda)})-\alpha\right\} \label{eq:k_formula1}\end{eqnarray}
and thus, since $z$ has generalised Cesaro limit $0$ and $\alpha\overset{C}{\sim}\frac{1}{2}C$
(by (1)(a)), so

\[
k\overset{C}{\sim}\frac{1}{C}\left\{ -i(s_{0}-r_{0}^{(\lambda)})-\frac{1}{2}C\right\} \]
as claimed in (1)(d).

Next, consider $k\alpha^{n}$. We have, in light of (1)(a), that

\begin{eqnarray*}
P[k(\alpha^{n}-\frac{1}{n+1}C^{n})] & = & \frac{1}{T}\left\{ \sum_{j=0}^{k-1}j\cdot0+k\left(\frac{\alpha^{n+1}}{n+1}-\frac{C^{n}}{n+1}\alpha\right)\right\} \\
\\ & = & \frac{1}{Ck+\tau_{0}+\alpha}\cdot k\left(\frac{\alpha^{n+1}}{n+1}-\frac{C^{n}}{n+1}\alpha\right)\\
\\ & = & \frac{1}{C}\left(\frac{\alpha^{n+1}}{n+1}-\frac{C^{n}}{n+1}\alpha\right)\left(1-\frac{(\tau_{0}+\alpha)}{Ck}+\ldots\right)\\
\\ & \overset{C}{\sim} & \frac{1}{C}\left(\frac{C^{n+1}}{(n+1)(n+2)}-\frac{1}{2}\frac{C^{n+1}}{n+1}\right)\\
\\ & = & -\frac{n}{2(n+1)(n+2)}C^{n}\end{eqnarray*}
It follows that

\begin{eqnarray}
k\alpha^{n} & \overset{C}{\sim} & \frac{C^{n}}{(n+1)}k-\frac{n}{2(n+1)(n+2)}C^{n}\nonumber \\
\nonumber \\ & \overset{C}{\sim} & \frac{C^{n-1}}{(n+1)}\left\{ -i(s_{0}-r_{0}^{(\lambda)})-[\frac{1}{2}C+\frac{n}{2(n+2)}C]\right\} \nonumber \\
\nonumber \\ & \overset{C}{\sim} & -i\,\frac{C^{n-1}}{(n+1)}(s_{0}-r_{0}^{(\lambda)})-\frac{1}{(n+2)}C^{n}\label{eq:k_alpha_n_formula1}\end{eqnarray}
using the result just derived for the Cesaro limit of $k$; and hence
for $n=1$ and $n=2$ we get the first two results claimed in (1)(b).

From (\ref{eq:k_alpha_n_formula1}) it then follows immediately that 

\begin{eqnarray}
z\alpha^{n} & = & ((s_{0}-\sigma_{0})-iT)\alpha^{n}=((s_{0}-\sigma_{0})-i(Ck+\tau_{0}+\alpha))\alpha^{n}\nonumber \\
\nonumber \\ & = & (s_{0}-r_{0}^{(\lambda)})\alpha^{n}-i\alpha^{n+1}-iCk\alpha^{n}\nonumber \\
\nonumber \\ & \overset{C}{\sim} & (s_{0}-r_{0}^{(\lambda)})\frac{C^{n}}{n+1}-i\frac{C^{n+1}}{n+2}-\frac{C^{n}}{n+1}(s_{0}-r_{0}^{(\lambda)})+i\frac{C^{n+1}}{n+2}=0\label{eq:z_alpha_n_formula1}\end{eqnarray}
also verifying the first two results of (1)(c).

And thus, noting (\ref{eq:k_formula1}) and invoking the results already
derived in (1)(a) and (1)(c) and the fact that the Cesaro eigenfunctions
$z^{2}$ and $z$ have generalised Cesaro limit $0$, we have

\begin{eqnarray*}
k^{2} & = & \frac{1}{C^{2}}\left\{ \begin{array}{cc}
-z^{2}+2(s_{0}-r_{0}^{(\lambda)})z-2iz\alpha\\
\\-(s_{0}-r_{0}^{(\lambda)})^{2}+2i(s_{0}-r_{0}^{(\lambda)})\alpha+\alpha^{2}\end{array}\right\} \\
\\ & \overset{C}{\sim} & \frac{1}{C^{2}}\left\{ -(s_{0}-r_{0}^{(\lambda)})^{2}+iC(s_{0}-r_{0}^{(\lambda)})+\frac{1}{3}C^{2}\right\} \end{eqnarray*}
which verifies the second result claimed in (1)(d).

Next we consider $k^{2}\alpha^{n}$. Along similar lines to our earlier
calculation for $k\alpha^{n}$ we have that

\begin{eqnarray*}
P[k^{2}(\alpha^{n}-\frac{1}{n+1}C^{n})] & = & \frac{1}{Ck+\tau_{0}+\alpha}\left\{ \sum_{j=0}^{k-1}j\cdot0+k^{2}\left(\frac{\alpha^{n+1}}{n+1}-\frac{C^{n}}{n+1}\alpha\right)\right\} \\
\\ & = & \frac{k}{C}\cdot\left(\frac{\alpha^{n+1}}{n+1}-\frac{C^{n}}{n+1}\alpha\right)\left(1-\frac{(\tau_{0}+\alpha)}{Ck}+\ldots\right)\\
\\ & = & \frac{1}{C}\left\{ \begin{array}{cc}
\frac{k\alpha^{n+1}}{n+1}-\frac{C^{n}}{n+1}k\alpha-\frac{(\tau_{0}\alpha^{n+1}+\alpha^{n+2})}{C(n+1)}\\
\\+\frac{C^{n-1}}{(n+1)}(\tau_{0}\alpha+\alpha^{2})+\ldots\end{array}\right\} \end{eqnarray*}
Using (\ref{eq:k_alpha_n_formula1}) and the result from (1)(a) it
follows that 

\begin{eqnarray*}
P[k^{2}(\alpha^{n}-\frac{1}{n+1}C^{n})] & \overset{C}{\sim} & \frac{1}{C}\left\{ \begin{array}{cc}
\frac{1}{(n+1)}\left[-i\frac{C^{n}}{(n+2)}(s_{0}-r_{0}^{(\lambda)})-\frac{C^{n+1}}{(n+3)}\right]\\
\\-\frac{C^{n}}{n+1}\left[-\frac{i}{2}(s_{0}-r_{0}^{(\lambda)})-\frac{C}{3}\right]\\
\\-\frac{1}{C(n+1)}\left[\tau_{0}\frac{C^{n+1}}{n+2}+\frac{C^{n+2}}{n+3}\right]\\
\\+\frac{C^{n-1}}{(n+1)}\left[\frac{\tau_{0}C}{2}+\frac{C^{2}}{3}\right]\end{array}\right\} \\
\\ & = & \frac{1}{C}\left\{ \begin{array}{cc}
iC^{n}\frac{n}{2(n+1)(n+2)}(s_{0}-r_{0}^{(\lambda)})\\
\\+\frac{2n}{3(n+1)(n+3)}C^{n+1}\\
\\+\frac{n}{2(n+1)(n+2)}\tau_{0}C^{n}\end{array}\right\} \end{eqnarray*}
and thus, by the result just proved for the Cesaro limit of $k^{2}$,
we have that

\begin{eqnarray}
k^{2}\alpha^{n} & \overset{C}{\sim} & \frac{C^{n}}{n+1}k^{2}+P[k^{2}(\alpha^{n}-\frac{1}{n+1}C^{n})]\nonumber \\
\nonumber \\ & \overset{C}{\sim} & \left\{ \begin{array}{cc}
-\frac{C^{n-2}}{(n+1)}(s_{0}-r_{0}^{(\lambda)})^{2}\\
\\+iC^{n-1}\frac{(3n+4)}{2(n+1)(n+2)}(s_{0}-r_{0}^{(\lambda)})\\
\\+\frac{C^{n}}{(n+3)}+\frac{n}{2(n+1)(n+2)}\tau_{0}C^{n-1}\end{array}\right\} \label{eq:k_2_alpha_n_formula1}\end{eqnarray}
which for $n=1$ proves the last of the claimed results in (1)(b).

Hence, along the same lines as our previous calculation for $z\alpha^{n}$,
we have 

\begin{eqnarray}
z^{2}\alpha^{n} & = & ((s_{0}-r_{0}^{(\lambda)})-iCk-i\alpha)^{2}\alpha^{n}\nonumber \\
\nonumber \\ & = & \left\{ \begin{array}{cc}
(s_{0}-r_{0}^{(\lambda)})^{2}\alpha^{n}-C^{2}k^{2}\alpha^{n}-\alpha^{n+2}\\
\\-2i(s_{0}-r_{0}^{(\lambda)})Ck\alpha^{n}-2Ck\alpha^{n+1}-2i(s_{0}-r_{0}^{(\lambda)})\alpha^{n+1}\end{array}\right\} \nonumber \\
\nonumber \\ & \overset{C}{\sim} & \frac{n}{2(n+1)(n+2)}C^{n+1}\left\{ i(s_{0}-r_{0}^{(\lambda)})-\tau_{0}\right\} \nonumber \\
\nonumber \\ & = & \frac{n}{2(n+1)(n+2)}C^{n+1}i(s_{0}-\sigma_{0})\label{eq:z_2_alpha_n_formula1}\end{eqnarray}
after invoking all the Cesaro limits already derived above and extensive
cancellations. When $n=1$ this verifies the last of the claimed results
in (1)(c).

This leaves only the Cesaro limit of $k^{3}$ to calculate. Here,
mimicking the calculation for $k^{2}$ and using the results from
(1)(a) and (1)(c) and the fact that the Cesaro eigenfunctions $z,\: z^{2}$
and $z^{3}$ all have generalised Cesaro limit $0$, we have that

\begin{eqnarray*}
k^{3} & = & \frac{1}{C^{3}}\left\{ \begin{array}{cc}
-iz^{3}+3z^{2}[i(s_{0}-r_{0}^{(\lambda)})+\alpha]\\
\\+3iz[i(s_{0}-r_{0}^{(\lambda)})+\alpha]^{2}\\
\\-[i(s_{0}-r_{0}^{(\lambda)})+\alpha]^{3}\end{array}\right\} \\
\\ & \overset{C}{\sim} & \frac{1}{C^{3}}\left\{ \begin{array}{cc}
\frac{i}{4}C^{2}(s_{0}-\sigma_{0})+i(s_{0}-r_{0}^{(\lambda)})^{3}+\frac{3C}{2}(s_{0}-r_{0}^{(\lambda)})^{2}\\
\\-iC^{2}(s_{0}-r_{0}^{(\lambda)})-\frac{1}{4}C^{3}\end{array}\right\} \end{eqnarray*}
which verifies the last of the claimed results in (1)(d) and thus
completes the proof.

Using this lemma we now calculate $r_{\zeta_{k}}^{(\lambda)}(s_{0},\mu)$
for $\mu=0,-1$ and $-2$ by Cesaro means. Let $N_{+,\mu}(T)$ be
the partial-sum function arising from summation over roots with imaginary
parts in $[0,T)$ on the line $\Re(s)=\sigma_{0}$ in the definition
of $r_{\zeta_{k}}^{(\lambda)}(s_{0},\mu)$; and let $N_{-,\mu}(T)$
be the corresponding partial-sum function for the roots with imaginary
parts in $(-T,0)$. Thus $N_{+,\mu}(T)$ is actually the partial-sum
function on the contour traversed by $z=(s_{0}-\sigma_{0})-iT$, given
by

\begin{equation}
N_{+,\mu}(T)=\sum_{j=0}^{k}((s_{0}-r_{0}^{(\lambda)})-iCj)^{-\mu}\label{eq:N_plus_mu_defn}\end{equation}
and $N_{-,\mu}(T)$ is the partial-sum function on the contour traversed
by $\tilde{z}=(s_{0}-\sigma_{0})+i\tilde{T}$, given by

\begin{equation}
N_{-,\mu}(\tilde{T})=\sum_{j=1}^{\tilde{k}}((s_{0}-r_{0}^{(\lambda)})+iCj)^{-\mu}\label{eq:N_minus_mu_defn}\end{equation}
It follows at once that, in general,

\begin{equation}
r_{\zeta_{k}}^{(\lambda)}(s_{0},\mu)=\textrm{e}^{i\pi\mu}\nu_{\lambda}\cdot\underset{z,\tilde{z}\rightarrow\infty}{Clim}\left\{ N_{+,\mu}(T)+N_{-,\mu}(\tilde{T})\right\} \label{eq:r_zeta_k_lambda_formula1}\end{equation}

\paragraph{Case 1 - $\mu=0$\textmd{:}}

In this case

\[
N_{+,0}(T)=k+1\overset{C}{\sim}-\frac{i}{C}(s_{0}-r_{0}^{(\lambda)})+\frac{1}{2}\]
and

\[
N_{-,0}(\tilde{T})=\tilde{k}\overset{C}{\sim}\frac{i}{C}(s_{0}-r_{0}^{(\lambda)})-\frac{1}{2}\]
Thus immediately in (\ref{eq:r_zeta_k_lambda_formula1}) we have

\begin{equation}
r_{\zeta_{k}}^{(\lambda)}(s_{0},0)=0\label{eq:r_zeta_k_lambda_s0_0}\end{equation}
as a function of $s_{0}$, and therefore in (\ref{eq:zeta_k_RootSide_Sum_lambdas}) 

\begin{equation}
r_{\zeta_{k}}(s_{0},0)=0\label{eq:r_zeta_k_s0_0}\end{equation}
also as a function of $s_{0}$, so that $d_{\zeta_{k}}(s_{0},0)=r_{\zeta_{k}}(s_{0},0)$
and the generalised root identity for $\zeta_{k}$ is confirmed when
$\mu=0$.

\paragraph{Case 2 - $\mu=-1$\textmd{:}}

Here, after simplification,

\begin{eqnarray*}
N_{+,-1}(T) & = & (s_{0}-r_{0}^{(\lambda)})(k+1)-i\frac{C}{2}(k^{2}+k)\\
\\ & \overset{C}{\sim} & -\frac{i}{2C}(s_{0}-r_{0}^{(\lambda)})^{2}+\frac{1}{2}(s_{0}-r_{0}^{(\lambda)})+\frac{i}{12}C\end{eqnarray*}
and

\begin{eqnarray*}
N_{-,-1}(\tilde{T}) & = & (s_{0}-r_{0}^{(\lambda)})\tilde{k}+i\frac{C}{2}(\tilde{k}^{2}+\tilde{k})\\
\\ & \overset{C}{\sim} & \frac{i}{2C}(s_{0}-r_{0}^{(\lambda)})^{2}-\frac{1}{2}(s_{0}-r_{0}^{(\lambda)})-\frac{i}{12}C\end{eqnarray*}
Thus again we have immediately in (\ref{eq:r_zeta_k_lambda_formula1})
that

\begin{equation}
r_{\zeta_{k}}^{(\lambda)}(s_{0},-1)=0\label{eq:r_zeta_k_lambda_s0_-1}\end{equation}
as a function of $s_{0}$, and therefore in (\ref{eq:zeta_k_RootSide_Sum_lambdas}) 

\begin{equation}
r_{\zeta_{k}}(s_{0},-1)=0\label{eq:r_zeta_k_s0_-1}\end{equation}
also as a function of $s_{0}$, so that likewise $d_{\zeta_{k}}(s_{0},-1)=r_{\zeta_{k}}(s_{0},-1)$
as expected and the generalised root identity for $\zeta_{k}$ is
confirmed when $\mu=-1$.

\paragraph{Case 3 - $\mu=-1$\textmd{:}}

In this case, after extensive simplification,

\begin{eqnarray*}
N_{+,-2}(T) & = & \left\{ \begin{array}{cc}
(s_{0}-r_{0}^{(\lambda)})^{2}(k+1)-iC(s_{0}-r_{0}^{(\lambda)})(k^{2}+k)\\
\\-C^{2}(\frac{1}{3}k^{3}+\frac{1}{2}k^{2}+\frac{1}{6}k)\end{array}\right\} \\
\\ & \overset{C}{\sim} & -\frac{i}{3C}(s_{0}-r_{0}^{(\lambda)})^{3}+\frac{1}{2}(s_{0}-r_{0}^{(\lambda)})^{2}+\frac{i}{12}C(s_{0}-r_{0}^{(\lambda)})+\frac{1}{12}C\tau_{0}\end{eqnarray*}
and

\begin{eqnarray*}
N_{-,-2}(\tilde{T}) & = & \left\{ \begin{array}{cc}
(s_{0}-r_{0}^{(\lambda)})^{2}\tilde{k}+iC(s_{0}-r_{0}^{(\lambda)})(\tilde{k}^{2}+\tilde{k})\\
\\-C^{2}(\frac{1}{3}\tilde{k}^{3}+\frac{1}{2}\tilde{k}^{2}+\frac{1}{6}\tilde{k})\end{array}\right\} \\
\\ & \overset{C}{\sim} & \frac{i}{3C}(s_{0}-r_{0}^{(\lambda)})^{3}-\frac{1}{2}(s_{0}-r_{0}^{(\lambda)})^{2}-\frac{i}{12}C(s_{0}-r_{0}^{(\lambda)})-\frac{1}{12}C\tau_{0}\end{eqnarray*}
Hence again in (\ref{eq:r_zeta_k_lambda_formula1}) we have immediately
that

\begin{equation}
r_{\zeta_{k}}^{(\lambda)}(s_{0},-2)=0\label{eq:r_zeta_k_lambda_s0_-2}\end{equation}
as a function of $s_{0}$, and therefore in (\ref{eq:zeta_k_RootSide_Sum_lambdas}) 

\begin{equation}
r_{\zeta_{k}}(s_{0},-2)=0\label{eq:r_zeta_k_s0_-2}\end{equation}
as a function of $s_{0}$, so that once more $d_{\zeta_{k}}(s_{0},-2)=r_{\zeta_{k}}(s_{0},-2)$
and the generalised root identity for $\zeta_{k}$ is also confirmed
when $\mu=-2$.

Overall then we confirm that, for the setting of zeta functions of
curves over finite fields, the root identities are satisfied when
$\mu=0,-1$ and $-2$. Moreover, it is clear from the general form
of such zeta functions in (\ref{eq:zeta_k_general_form2}) and the
results derived above for $r_{\zeta_{k}}^{(\lambda)}(s_{0},\mu)$
when $\mu=0,-1$ and $-2$ that in this setting the generalised root
identities could not be used to detect any possible violation of the
RH for such zeta functions, since these results for $r_{\zeta_{k}}^{(\lambda)}$
apply for arbitrary $\lambda$ with $\sigma_{0}$ any value in $[0,1]$.

Of course in fact the RH is known to hold for such zeta functions
of curves over finite fields, but the above observation nonetheless
demonstrates one way of seeing why this state of affairs is consistent
with the generalised root identities, in contrast to the situation
for the Riemann zeta function discussed in {[}1{]}, where the root
identity for $\mu=-2$ is argued to imply a contradiction of the RH.

In the next section, we provide a second way of seeing this that more
closely identifies this contrast. We re-perform the calculations for
$r_{\zeta_{k}}(s_{0},\mu)$ when $\mu=0,-1$ and $-2$ in a way that
directly mimics the calculation for $\zeta$ in {[}1{]}, by working
solely on the line $\Re(s)=\frac{1}{2}$ and considering all the roots
simultaneously rather than splitting them into a union of sets of
equally spaced roots as just done.

\section{A second approach to the Cesaro calculation of $r_{\zeta_{k}}(s_{0},\mu)$
for $\mu=0,-1$ and $-2$}

For simplicity we shall use the results for $r_{\zeta_{k}}^{(\lambda)}(s_{0},\mu)$
($\mu=0,-1,-2$) from the last section just to allow us to ignore
the contributions to $r_{\zeta_{k}}(s_{0},\mu)$ from the poles which
lie equally spaced on $\Re(s)=0$ and $\Re(s)=1$, and thus to concentrate
solely on the roots (arising from the factorisation of $P_{2g}(q^{-s})$)
to which the RH applies.

We let $N(T)$ be the counting function giving the number of roots
of $P_{2g}(q^{-s})$ in the critical strip with imaginary part in
$(0,T)$. Since the roots of $\zeta_{k}$ are mirror-symmetric in
the real axis, so the number of roots in the critical strip with imaginary
part in $(-\tilde{T},0)$ is likewise $N(\tilde{T})$. Here, without
loss of generality, we are assuming that there are no base roots with
imaginary part precisely $0$ - if there were they could be ignored,
just as we have ignored the poles, as constituting a further discrete
set of equi-spaced roots of the sort analyzed in section 3.%
\footnote{Alternatively we could readily incorporate this case into the analysis
in this section, as we could with the poles, but for simplicity we
opt not to do this here%
}

In this case we add $2g$ roots every time $T$ increases by $C$,
so

\begin{equation}
N(T)=\frac{2g}{C}T+S(T)\label{eq:N_T_defn1}\end{equation}
where $S(T)$ is then the periodic, period-$C$ function obtained
by periodically extending the function $N(T)-\frac{2g}{C}T$ on $[0,C)$
to all of $[0,\infty)$. In analogy with the notation used in {[}1{]}
for $\zeta$ we denote the divergent piece of $N(T)$ by

\begin{equation}
\check{N}(T)=\frac{2g}{C}T\label{eq:N_check_T_defn1}\end{equation}
and $S(T)$ plays the role analogous to the oscillatory argument of
the Riemann zeta function. In the case of $\zeta_{k}$ there is no
asymptotically decaying piece analogous to the term $\frac{1}{\pi}\delta(T)$.

Note that $S(0)=S(C)=0$ and that, since the conjugate of a root of
$\zeta_{k}$ is also a root, so $S(t)=-S(C-t)$ on $[0,C)$. Thus
$S$ is an odd function about the centre-point of each period and
pinned at $0$ at the end points.

This information, applied in the breakdown given in (\ref{eq:N_T_defn1}),
is already sufficient to calculate $r_{\zeta_{k}}(s_{0},\mu)$ for
$\mu=0$ and $-1$, as follows. As in {[}1{]}, we treat the functions
$N(T)$ and $N(\tilde{T})$ as being functions on the critical line
(either as a consequence of the RH for $\zeta_{k}$ or, if this is
not to be assumed here, as a reflection of the mirror-symmetry of
roots of $\zeta_{k}$ in the critical line); thus we have

\begin{equation}
z=(s_{0}-\frac{1}{2})-iT\qquad\textrm{so}\qquad T=i(z-(s_{0}-\frac{1}{2}))\label{eq:z_defn_critical_line}\end{equation}
and

\begin{equation}
\tilde{z}=(s_{0}-\frac{1}{2})+i\tilde{T}\qquad\textrm{so}\qquad\tilde{T}=-i(\tilde{z}-(s_{0}-\frac{1}{2}))\label{eq:ztilde_defn_critical_line}\end{equation}

\paragraph{Case 1 - $\mu=0$:}

In this case, since all roots have $M_{i}=1$, we have

\begin{eqnarray}
r_{\zeta_{k}}(s_{0},0) & = & \sum_{\{s_{0}-roots\:\rho_{i}\: of\:\zeta_{k}\}}(s_{0}-\rho_{i})^{0}=\underset{z,\tilde{z}\rightarrow\infty}{Clim}\,\{N(T)+N(\tilde{T})\}\nonumber \\
\nonumber \\ & = & \underset{z,\tilde{z}\rightarrow\infty}{Clim}\,\{[\check{N}(T)+S(T)]+[\check{N}(\tilde{T})+S(\tilde{T})]\}\label{eq:r_zeta_k_s0_0_Alt}\end{eqnarray}
Now, since $S$ is an odd function about its mid-point on each period,
so $\int_{0}^{C}S(t)\,\textrm{d}t=0$ and its average value on each
period is zero. It follows immediately that

\[
\underset{T\rightarrow\infty}{Clim}\, S(T)=\underset{\tilde{T}\rightarrow\infty}{Clim}\, S(\tilde{T})=0\]
And since $\underset{z\rightarrow\infty}{Clim}\, z=\underset{\tilde{z}\rightarrow\infty}{Clim}\,\tilde{z}=0$,
so 

\[
\underset{z\rightarrow\infty}{Clim}\,\check{N}(T)=\underset{z\rightarrow\infty}{Clim}\,\frac{2g}{C}i(z-(s_{0}-\frac{1}{2}))=-\frac{2g}{C}i(s_{0}-\frac{1}{2})\]
and

\[
\underset{\tilde{z}\rightarrow\infty}{Clim}\,\check{N}(\tilde{T})=\underset{\tilde{z}\rightarrow\infty}{Clim}\,-\frac{2g}{C}i(\tilde{z}-(s_{0}-\frac{1}{2}))=\frac{2g}{C}i(s_{0}-\frac{1}{2})\]
Combining these results in (\ref{eq:r_zeta_k_s0_0_Alt}) it follows
at once that 

\[
r_{\zeta_{k}}(s_{0},0)=0\]
as claimed before in Section 3 and in agreement with the root identity
for $\zeta_{k}$ at $\mu=0$ and Result 2 for the value of $d_{\zeta_{k}}(s_{0},0)$.

\paragraph{Case 2 - $\mu=-1$:}

In this case, since the roots $\rho_{i}=\beta_{i}+i\gamma_{i}$ of
$\zeta_{k}$ are mirror-symmetric in the critical line, we have, on
invoking the result just proved for $r_{\zeta_{k}}(s_{0},0)$, that

\begin{eqnarray*}
r_{\zeta_{k}}(s_{0},-1) & = & -\sum_{\{s_{0}-roots\:\rho_{i}\: of\:\zeta_{k}\}}(s_{0}-\rho_{i})^{1}\\
\\ & = & -(s_{0}-\frac{1}{2})r_{\zeta_{k}}(s_{0},0)+i\cdot\sum_{\{s_{0}-roots\:\rho_{i}\: of\:\zeta_{k}\}}\gamma_{i}\\
\\ & = & i\cdot\underset{z,\tilde{z}\rightarrow\infty}{Clim}\,\left\{ \int_{0}^{T}t\,\textrm{d}N(t)-\int_{0}^{\tilde{T}}\tilde{t}\,\textrm{d}N(\tilde{t})\right\} \end{eqnarray*}
Using the notation from {[}1{]} under which $N_{i}(T)$ represents
the i'th integral of $N(T)$ and similarly for $\check{N}(T)$ and
$S(T)$, this becomes

\begin{eqnarray}
 &  & r_{\zeta_{k}}(s_{0},-1)\nonumber \\
\nonumber \\ & = & i\cdot\underset{z,\tilde{z}\rightarrow\infty}{Clim}\,\left\{ \left\{ TN(T)-N_{1}(T)\right\} -\left\{ \tilde{T}N(\tilde{T})-N_{1}(\tilde{T})\right\} \right\} \nonumber \\
\nonumber \\ & = & i\cdot\underset{z,\tilde{z}\rightarrow\infty}{Clim}\,\left\{ \left\{ \begin{array}{cc}
\left\{ T\check{N}(T)-\check{N}_{1}(T)\right\} \\
\\+\left\{ TS(T)-S_{1}(T)\right\} \end{array}\right\} -\left\{ \begin{array}{cc}
\left\{ \tilde{T}\check{N}(\tilde{T})-\check{N}_{1}(\tilde{T})\right\} \\
\\+\left\{ \tilde{T}S(\tilde{T})-S_{1}(\tilde{T})\right\} \end{array}\right\} \right\} \nonumber \\
\label{eq:r_zeta_k_s0_-1_Alt}\end{eqnarray}
Now, as before,

\[
T\check{N}(T)=\frac{2g}{C}T^{2}\overset{C}{\rightarrow}-\frac{2g}{C}(s_{0}-\frac{1}{2})^{2}\]
and

\[
\check{N}_{1}(T)=\frac{g}{C}T^{2}\overset{C}{\rightarrow}-\frac{g}{C}(s_{0}-\frac{1}{2})^{2}\]
so

\[
T\check{N}(T)-\check{N}_{1}(T)\overset{C}{\rightarrow}-\frac{g}{C}(s_{0}-\frac{1}{2})^{2}\]
and in identical fashion

\[
\tilde{T}\check{N}(\tilde{T})-\check{N}_{1}(\tilde{T})\overset{C}{\rightarrow}-\frac{g}{C}(s_{0}-\frac{1}{2})^{2}.\]

As for the terms involving $S(T)$, since $\int_{0}^{C}S(t)\,\textrm{d}t=0$,
so $S_{1}(T)$ is also a periodic function with period $C$ satisfying
$S_{1}(0)=S_{1}(C)=0$. As such it converges in a generalised Cesaro
sense, under a single application of the averaging operator $P$,
to its average value on each period. Unlike $S(T)$, however, since
it is not odd on each period, this average value, $S_{1}^{av}$, is
not zero; i.e. overall

\[
\underset{T\rightarrow\infty}{Clim}\, S_{1}(T)=\underset{\tilde{T}\rightarrow\infty}{Clim}\, S_{1}(\tilde{T})=S_{1}^{av}\]
where

\[
S_{1}^{av}=\frac{1}{C}\int_{0}^{C}S_{1}(t)\,\textrm{d}t\neq0\]

This leaves only the terms $TS(T)$ and $\tilde{T}S(\tilde{T})$.
These both have generalised Cesaro limit zero, i.e. 

\[
\underset{T\rightarrow\infty}{Clim}\, TS(T)=\underset{\tilde{T}\rightarrow\infty}{Clim}\,\tilde{T}S(\tilde{T})=0\]
To see this we again simply apply $P$ and use integration by parts,
but taking care in doing so to use the unique anti-derivative of $S(t)$
having average value $0$, namely $S_{1}(t)-S_{1}^{av}$. Specifically

\begin{eqnarray*}
P[tS(t)](T) & = & \frac{1}{T}\int_{0}^{T}tS(t)\,\textrm{d}t\\
\\ & = & (S_{1}(T)-S_{1}^{av})-P[S_{1}(t)-S_{1}^{av}](T)\end{eqnarray*}
Since clearly $P[S_{1}(t)-S_{1}^{av}](T)\rightarrow0$ as $T\rightarrow\infty$,
so $P^{2}[tS(t)](T)\rightarrow0$ as $T\rightarrow\infty$ and thus
$\underset{T\rightarrow\infty}{Clim}\, TS(T)=0$ as claimed; and likewise
for $\underset{\tilde{T}\rightarrow\infty}{Clim}\,\tilde{T}S(\tilde{T})$.

Combining these calculations in (\ref{eq:r_zeta_k_s0_-1_Alt}), we
finally deduce that 

\begin{eqnarray*}
r_{\zeta_{k}}(s_{0},-1) & = & 0\end{eqnarray*}
once again, as claimed before in section 3 and in line with the root
identity for $\zeta_{k}$ at $\mu=-1$ and Result 2 for the value
of $d_{\zeta_{k}}(s_{0},-1)$.

\paragraph{Case 3 - $\mu=-2$:}

In this case, writing $\rho_{i}=\beta_{i}+i\gamma_{i}=(\frac{1}{2}+\epsilon_{i})+i\gamma_{i}$
as in {[}1{]}, and again invoking the mirror symmetry of roots of
$\zeta_{k}$ in the critical line and the results just derived for
$r_{\zeta_{k}}(s_{0},0)$ and $r_{\zeta_{k}}(s_{0},-1)$, we obtain

\begin{eqnarray}
r_{\zeta_{k}}(s_{0},-2) & = & \sum_{\{s_{0}-roots\:\rho_{i}\: of\:\zeta_{k}\}}(s_{0}-\rho_{i})^{2}\nonumber \\
\nonumber \\ & = & \sum_{\{s_{0}-roots\:\rho_{i}\: of\:\zeta_{k}\}}((s_{0}-\frac{1}{2})-\epsilon_{i}-i\gamma_{i})^{2}\nonumber \\
\nonumber \\ & = & \left\{ \begin{array}{cc}
(s_{0}-\frac{1}{2})^{2}r_{\zeta_{k}}(s_{0},0)+X_{\epsilon}-2(s_{0}-\frac{1}{2})r_{\zeta_{k}}(s_{0},-1)\\
\\-\sum_{\{s_{0}-roots\:\rho_{i}\: of\:\zeta_{k}\}}\gamma_{i}^{2}\end{array}\right\} \nonumber \\
\nonumber \\ & = & X_{\epsilon}-\underset{z,\tilde{z}\rightarrow\infty}{Clim}\,\left\{ \int_{0}^{T}t^{2}\,\textrm{d}N(t)+\int_{0}^{\tilde{T}}\tilde{t^{2}}\,\textrm{d}N(\tilde{t})\right\} \label{eq:r_zeta_k_s0_-2_AltA}\end{eqnarray}
where $X_{\epsilon}=\sum_{\{s_{0}-roots\:\rho_{i}\: of\:\zeta_{k}\}}\epsilon_{i}^{2}$
as per the notation of {[}1{]} (and thus $X_{\epsilon}$ is trivially
zero under the RH for $\zeta_{k}$).

Thus, on integration by parts,

\begin{eqnarray}
 &  & r_{\zeta_{k}}(s_{0},-2)\nonumber \\
\nonumber \\ & = & X_{\epsilon}-\underset{z,\tilde{z}\rightarrow\infty}{Clim}\,\left\{ \begin{array}{cc}
\left\{ T^{2}N(T)-2TN_{1}(T)+2N_{2}(T)\right\} \\
\\+\left\{ \tilde{T}^{2}N(\tilde{T})-2\tilde{T}N_{1}(\tilde{T})+2N_{2}(\tilde{T})\right\} \end{array}\right\} \nonumber \\
\nonumber \\ & = & X_{\epsilon}-\underset{z,\tilde{z}\rightarrow\infty}{Clim}\,\left\{ \begin{array}{cc}
\left\{ T^{2}\check{N}(T)-2T\check{N}_{1}(T)+2\check{N}_{2}(T)\right\} \\
\\+\left\{ T^{2}S(T)-2TS_{1}(T)+2S_{2}(T)\right\} \\
\\+\left\{ \tilde{T}^{2}\check{N}(\tilde{T})-2\tilde{T}\check{N}_{1}(\tilde{T})+2\check{N}_{2}(\tilde{T})\right\} \\
\\+\left\{ \tilde{T}^{2}S(\tilde{T})-2\tilde{T}S_{1}(\tilde{T})+2S_{2}(\tilde{T})\right\} \end{array}\right\} \label{eq:r_zeta_k_s0_-2_AltB}\end{eqnarray}
For this case, however, unlike for $\mu=0$ and $\mu=-1$, we are
not yet in a position to complete the calculation of $r_{\zeta_{k}}(s_{0},-2)$,
in particular because for this case we need to calculate $S_{1}^{av}$
(and then calculate for $S_{2}$).

To proceed, therefore, we note that it suffices simply to consider
the case in which $g=1$ and there are precisely two distinct base
roots with imaginary parts at $\kappa$ and $C-\kappa$ in the initial
period $0<T<C$ on the critical line. This suffices because the case
of general $N(T)$ arising from $2g$ base roots in $0\leq T<C$ (either
with $2g$ distinct steps of height $1$ or possibly fewer steps,
some of height greater than $1$ arising from base roots with the
same imaginary part) can clearly be made up as a linear combination
of such simpler cases, and since the expression (\ref{eq:r_zeta_k_s0_-2_AltA})
is clearly linear in $N(T)$ and $N(\tilde{T})$, so the value of
$r_{\zeta_{k}}(s_{0},-2)$ in the general case is just the same linear
combination of the values of $r_{\zeta_{k}}(s_{0},-2)$ for these
simpler cases.%
\footnote{Note also that the case of $g=1$, but with $N(T)$ having a single
step of height $2$ at $\frac{C}{2}$, follows as a limiting case
of the one we consider in which $\kappa\rightarrow\frac{C}{2}$; so
in treating only this single simple case with distinct jumps at $\kappa$
and $C-\kappa$ we are indeed calculating all that is required in
order to cover the general case.%
}

We thus now assume without loss of generality that $g=1$ and $N(T)$
has precisely two distinct steps of height $1$, at $\kappa$ and
$C-\kappa$ on the interval $0\leq T<C$ on the critical line, i.e.
$N(t)$ has the form

\begin{eqnarray}
N(t) & = & \begin{cases}
0 & ,\:0\leq t<\kappa\\
1 & ,\:\kappa\leq t<C-\kappa\\
2 & ,\: C-\kappa\leq t\leq C\end{cases}\label{eq:N_T_defnA}\end{eqnarray}
on $0\leq t\leq C$. Thus, writing $T=Ck+\alpha$, $0\leq\alpha<C$
(and $\tilde{T}=C\tilde{k}+\tilde{\alpha}$) in the usual fashion,
we have

\begin{equation}
\check{N}(T)=\frac{2}{C}T\label{eq:N_check_T_defnA}\end{equation}
and $S(T)$ is the periodic, period-$C$ function with the properties
described earlier given by

\begin{eqnarray}
S(T) & = & -\frac{2}{C}\alpha+\begin{cases}
0 & ,\:0\leq\alpha<\kappa\\
1 & ,\:\kappa\leq\alpha<C-\kappa\\
2 & ,\: C-\kappa\leq\alpha\leq C\end{cases}\label{eq:S_T_defnA}\end{eqnarray}
Hence $S_{1}(T)$ is given by

\begin{eqnarray}
S_{1}(T) & = & -\frac{1}{C}\alpha^{2}+\begin{cases}
0 & ,\:0\leq\alpha<\kappa\\
\alpha-\kappa & ,\:\kappa\leq\alpha<C-\kappa\\
2\alpha-C & ,\: C-\kappa\leq\alpha\leq C\end{cases}\label{eq:S_1_T_defnA}\end{eqnarray}
and an elementary computation then shows that in this case

\begin{equation}
S_{1}^{av}=\frac{1}{C}\int_{0}^{C}S_{1}(t)\,\textrm{d}t=\frac{1}{C}\kappa^{2}-\kappa+\frac{1}{6}C\label{eq:S_1_av_calculn}\end{equation}

Taking the terms in (\ref{eq:r_zeta_k_s0_-2_AltB}) in turn then we
have, in the usual way, that

\begin{eqnarray*}
T^{2}\check{N}(T) & = & \frac{2}{C}T^{3}\overset{C}{\rightarrow}\frac{2}{C}i(s_{0}-\frac{1}{2})^{3}\\
\\T\check{N_{1}}(T) & = & \frac{1}{C}T^{3}\overset{C}{\rightarrow}\frac{1}{C}i(s_{0}-\frac{1}{2})^{3}\\
\\\check{N_{2}}(T) & = & \frac{1}{3C}T^{3}\overset{C}{\rightarrow}\frac{1}{3C}i(s_{0}-\frac{1}{2})^{3}\end{eqnarray*}
and similarly

\begin{eqnarray*}
\tilde{T}^{2}\check{N}(\tilde{T}) & \overset{C}{\rightarrow} & -\frac{2}{C}i(s_{0}-\frac{1}{2})^{3}\\
\\\tilde{T}\check{N_{1}}(\tilde{T}) & \overset{C}{\rightarrow} & -\frac{1}{C}i(s_{0}-\frac{1}{2})^{3}\\
\\\check{N_{2}}(\tilde{T}) & \overset{C}{\rightarrow} & -\frac{1}{3C}i(s_{0}-\frac{1}{2})^{3}\end{eqnarray*}
Thus

\begin{equation}
\underset{z,\tilde{z}\rightarrow\infty}{Clim}\,\left\{ \begin{array}{cc}
\left\{ T^{2}\check{N}(T)-2T\check{N}_{1}(T)+2\check{N}_{2}(T)\right\} \\
\\+\left\{ \tilde{T}^{2}\check{N}(\tilde{T})-2\tilde{T}\check{N}_{1}(\tilde{T})+2\check{N}_{2}(\tilde{T})\right\} \end{array}\right\} =0\label{eq:r_zeta_k_s0_-2_Calc_Nterms}\end{equation}
As for the terms in (\ref{eq:r_zeta_k_s0_-2_AltB}) involving $S(T)$,
since $S_{1}^{av}\neq0$, so $S_{2}(T)$ is no longer periodic and
thus requires the removal of a divergent piece before averaging in
order to find its Cesaro limit. Since we still have $\alpha^{n}\overset{C}{\sim}\frac{C^{n}}{n+1}$
as in section 3, and thus also

\[
Ck=T-\alpha\overset{C}{\sim}-i(s_{0}-\frac{1}{2})-\frac{1}{2}C\]
and likewise

\[
C\tilde{k}=\tilde{T}-\tilde{\alpha}\overset{C}{\sim}i(s_{0}-\frac{1}{2})-\frac{1}{2}C\]
so integrating the expression for $S_{1}(T)$ we have

\begin{eqnarray}
S_{2}(T) & = & S_{1}^{av}\cdot Ck+\int_{0}^{\alpha}S_{1}(t)\,\textrm{d}t\nonumber \\
\nonumber \\ & = & S_{1}^{av}\cdot Ck-\frac{1}{3C}\alpha^{3}+Q(\alpha)\nonumber \\
\nonumber \\ & \overset{C}{\sim} & \{-i(s_{0}-\frac{1}{2})-\frac{1}{2}C\}S_{1}^{av}-\frac{1}{12}C^{2}+Q(\alpha)\label{eq:S_2_InitCalc}\end{eqnarray}
where $Q$ is the (discontinuous), periodic, period-$C$ function
given by\[
Q(T)=Q(\alpha)=\begin{cases}
0 & ,\:0\leq\alpha<\kappa\\
\frac{1}{2}(\alpha-\kappa)^{2} & ,\:\kappa\leq\alpha<C-\kappa\\
(\alpha-\frac{1}{2}C)^{2}+(\kappa-\frac{1}{2}C)^{2} & ,\: C-\kappa\leq\alpha\leq C\end{cases}\]
Since $Q$ is periodic, its Cesaro limit is its average value, which
in light of (\ref{eq:S_1_av_calculn}) is given by

\begin{eqnarray}
Q^{av} & = & \frac{1}{C}\int_{0}^{C}Q(\alpha)\,\textrm{d}\alpha\nonumber \\
\nonumber \\ & = & \frac{1}{C}\left\{ \frac{1}{6}(C-2\kappa)^{3}+\frac{1}{3}\left\{ \frac{1}{8}C^{3}-(\frac{1}{2}C-\kappa)^{3}\right\} +(\kappa-\frac{1}{2}C)^{2}\kappa\right\} \nonumber \\
\nonumber \\ & = & \frac{1}{6}C^{2}-\frac{1}{2}C\kappa+\frac{1}{2}\kappa^{2}=\frac{1}{2}CS_{1}^{av}+\frac{1}{12}C^{2}\label{eq:Q_av_calculn}\end{eqnarray}
Combining (\ref{eq:Q_av_calculn}) in (\ref{eq:S_2_InitCalc}) we
thus have

\[
\underset{z\rightarrow\infty}{Clim}\, S_{2}(T)=-i(s_{0}-\frac{1}{2})S_{1}^{av}=-i(s_{0}-\frac{1}{2})\left\{ \frac{1}{C}\kappa^{2}-\kappa+\frac{1}{6}C\right\} \]
and likewise, by an identical computation,

\[
\underset{\tilde{z}\rightarrow\infty}{Clim}\, S_{2}(\tilde{T})=i(s_{0}-\frac{1}{2})S_{1}^{av}=i(s_{0}-\frac{1}{2})\left\{ \frac{1}{C}\kappa^{2}-\kappa+\frac{1}{6}C\right\} \]

We next calculate the Cesaro limit of $TS_{1}(T)$ (and its counterpart
$\tilde{T}S_{1}(\tilde{T})$); writing this as

\begin{equation}
TS_{1}(T)=S_{1}^{av}T+T(S_{1}(T)-S_{1}^{av})\label{eq:T_S_1_InitCalc}\end{equation}
we have, in the usual way,

\begin{equation}
S_{1}^{av}T\overset{C}{\rightarrow}-i(s_{0}-\frac{1}{2})S_{1}^{av}\label{eq:S_1_av_T_calc}\end{equation}
and it remains only to calculate the Cesaro limit of the second term
in (\ref{eq:T_S_1_InitCalc}). Since $S_{1}(T)-S_{1}^{av}$ is periodic
with average value $0$ its Cesaro limit is obtained from application
of a pure power of $P$, without need for removal of any further divergences.
Defining $\breve{S}_{2}(T):=\int_{0}^{T}(S_{1}(t)-S_{1}^{av})\,\textrm{d}t=S_{2}(T)-S_{1}^{av}T$,
on integration by parts we have

\begin{eqnarray}
P[t(S_{1}(t)-S_{1}^{av})](T) & = & \frac{1}{T}\int_{0}^{T}t(S_{1}(t)-S_{1}^{av})\,\textrm{d}t\nonumber \\
\nonumber \\ & = & \breve{S}_{2}(T)-P[\breve{S}_{2}](T)\label{eq:S_2_upper_calc}\end{eqnarray}
But the fact that $S_{1}(T)-S_{1}^{av}$ is periodic with average
value $0$ means that $\breve{S}_{2}(T)$ must also be periodic, with
Cesaro limit given by its average value in each period, $\breve{S}_{2}^{av}$,
i.e. 

\[
P[\breve{S}_{2}](T)\rightarrow\breve{S}_{2}^{av}\qquad\textrm{as}\qquad T\rightarrow\infty\]
It follows immediately in (\ref{eq:S_2_upper_calc}) that

\begin{equation}
P^{2}[t(S_{1}(t)-S_{1}^{av})](T)\rightarrow0\qquad\textrm{as}\qquad T\rightarrow\infty\label{eq:T_S_1_InterimCalcA}\end{equation}
and so

\begin{equation}
T(S_{1}(T)-S_{1}^{av})\overset{C}{\sim}0\label{eq:T_S_1_InterimCalcB}\end{equation}
Combining (\ref{eq:S_1_av_T_calc}) and (\ref{eq:T_S_1_InterimCalcB})
in (\ref{eq:T_S_1_InitCalc}) it follows that 

\[
\underset{z\rightarrow\infty}{Clim}\, TS_{1}(T)=-i(s_{0}-\frac{1}{2})S_{1}^{av}\]
and likewise, by an identical computation,

\[
\underset{\tilde{z}\rightarrow\infty}{Clim}\,\tilde{T}S_{1}(\tilde{T})=i(s_{0}-\frac{1}{2})S_{1}^{av}\]

We now lastly calculate the Cesaro limit of $T^{2}S(T)$ (and its
counterpart $\tilde{T}^{2}S(\tilde{T})$). Using $S_{1}(T)-S_{1}^{av}$
as an antiderivative of $S(T)$ in integration by parts, we have

\begin{eqnarray*}
P[t^{2}S(t)](T) & = & \frac{1}{T}\int_{0}^{T}t^{2}S(t)\,\textrm{d}t\\
\\ & = & T(S_{1}(T)-S_{1}^{av})-2P[t(S_{1}(t)-S_{1}^{av})](T)\end{eqnarray*}
But then it follows directly from (\ref{eq:T_S_1_InterimCalcA}) that

\[
P^{3}[t^{2}S(t)](T)\rightarrow0\qquad\textrm{as}\qquad T\rightarrow\infty\]
and so

\[
\underset{T\rightarrow\infty}{Clim}\, T^{2}S(T)=0\]
without any need to remove divergences prior to averaging.%
\footnote{Hence the use of the notation $\underset{T\rightarrow\infty}{Clim}$
rather than $\underset{z\rightarrow\infty}{Clim}$, as we have done
at times earlier. The lack of need to remove divergences reflects
the fact that $S(T)$ is odd about its mid-point on each period, so
that $T^{2}S(T)$ is essentially pure oscillatory, albeit of increasing
amplitude as $T\rightarrow\infty$%
}By identical reasoning we likewise have

\[
\underset{\tilde{T}\rightarrow\infty}{Clim}\,\tilde{T}^{2}S(\tilde{T})=0\]

Finally, combining all our results for $S(T)$-related terms in (\ref{eq:r_zeta_k_s0_-2_AltB}),
we have

\begin{equation}
\underset{z,\tilde{z}\rightarrow\infty}{Clim}\,\left\{ \begin{array}{cc}
\left\{ T^{2}S(T)-2TS_{1}(T)+2S_{2}(T)\right\} \\
\\+\left\{ \tilde{T}^{2}S(\tilde{T})-2\tilde{T}S_{1}(\tilde{T})+2S_{2}(\tilde{T})\right\} \end{array}\right\} =0\label{eq:r_zeta_k_s0_-2_Calc_Sterms}\end{equation}
(in fact each of the two pieces in this expression is itself zero)
and combining (\ref{eq:r_zeta_k_s0_-2_Calc_Nterms}) and (\ref{eq:r_zeta_k_s0_-2_Calc_Sterms})
in (\ref{eq:r_zeta_k_s0_-2_AltB}) it follows at last that

\begin{equation}
r_{\zeta_{k}}(s_{0},-2)=X_{\epsilon}\label{eq:r_zeta_k_s0_-2_AltC}\end{equation}

It follows at once from Result 2 for $d_{\zeta_{k}}(s_{0},-2)$ and
the root identity for $\zeta_{k}$ at $\mu=-2$ that we must have

\begin{equation}
X_{\epsilon}=0\label{eq:X_epsilon_calc}\end{equation}

Thus we see at once that, unlike for $\zeta$ in {[}1{]}, there is
no contradiction between the $\mu=-2$ root identity and the RH for
$\zeta_{k}$, since under the RH for $\zeta_{k}$ all $\epsilon_{i}=0$
and so $X_{\epsilon}$ is trivially zero. This is in contrast to the
case for $\zeta$, where we had $d_{\zeta}(s_{0},-2)=0$ but $r_{\zeta}(s_{0},-2)=-\frac{1}{2}+X_{\epsilon}$,
so that the root identity for $\zeta$ at $\mu=-2$ implies $X_{\epsilon}=\frac{1}{2}\neq0$
and thus we cannot have $\epsilon_{i}=0$ for all roots of $\zeta$
as required by the RH for $\zeta$.

Note, however, that the fact that the $\mu=-2$ root identity implies
$X_{\epsilon}=0$ in the case of $\zeta_{k}$ is not sufficient alone
to actually \textit{deduce} the RH for such zeta functions of curves
over finite fields. This is because, as we have shown in section 3,
it is possible to have infinitely many roots off the critical line,
but still have the Cesaro limit defining $X_{\epsilon}$ for these
roots be zero - for example if all the roots of this form lie equi-spaced
in a vertical line with $\Re(s)=\sigma_{0}\neq\frac{1}{2}$, or comprise
a union of such sets of roots.

Thus deduction of the RH for zeta functions of curves over finite
fields requires more than just the root identities for $\mu=0,-1$
and $-2$, and indeed it requires work reaching down to algebraic
considerations concerning the field $k$ of a sort which have not
arisen in the root identity computations here.

\section{Summary and Observations}

Overall we have shown that for zeta functions of curves over finite
fields, $\zeta_{k}$, the generalised root identities (\ref{eq:GenRtId1})
are satisfied for arbitrary $s_{0}$ ($\Re(s_{0})>1$) and $\mu$
and we have calculated the detailed implications of this for $\mu=0,-1$
and $-2$ in two independent ways.

The first approach, in section 3, uses the fact that the generalised
roots of such $\zeta_{k}$ form equi-spaced sets on vertical lines
$\Re(s)=\sigma_{0},\:0\leq\sigma_{0}\leq1$, to deduce directly that
we do have $d_{\zeta_{k}}(s_{0},\mu)=r_{\zeta_{k}}(s_{0},\mu)=0$
for $\mu=0,-1$ and $-2$.%
\footnote{And in fact this will continue to be true for $\mu=-3,-4,\ldots$
since these sort of bidirectional equally-spaced remainder sums generally
give $0$ as their Cesaro sums for elementary, entire sum-functions%
}

The second approach in section 4 directly mimics the calculational
approach used for the Riemann zeta function in {[}1{]} and shows that
we have $d_{\zeta_{k}}(s_{0},\mu)=r_{\zeta_{k}}(s_{0},\mu)=0$ when
$\mu=0$ or $-1$, while $d_{\zeta_{k}}(s_{0},-2)=0$ and $r_{\zeta_{k}}(s_{0},-2)=X_{\epsilon}$,
so that we must have $X_{\epsilon}=0$. This is certainly consistent
with the RH for $\zeta_{k}$, which is famously true, in contrast
to the situation for $\zeta$ in {[}1{]} where the $\mu=-2$ root
identity is seen to imply a contradiction of the RH. As noted in section
4, however, the fact that $X_{\epsilon}=0$ is not sufficient to actually
deduce the RH for $\zeta_{k}$ on its own.

We conclude with a brief discussion of why the implications of the
generalised root identities for $\zeta_{k}$ and $\zeta$ differ (the
RH is true for $\zeta_{k}$ but appears false for $\zeta$ on the
basis of the $\mu=-2$ root identity) even though the zeta functions
in both settings satisfy both Euler and Hadamard product formulas.

The reason appears to be that in the case of zeta functions of curves
over finite fields the detailed structure of the counting function
$N(T)$ counting non-trivial roots in the critical strip is much simpler
than it is for $\zeta$. Specifically:

\textbf{(a)} For $\zeta_{k}$, $N(T)$ has only divergent and oscillatory
pieces, $\check{N}(T)$ and $S(T)$, but no decaying asymptotic piece
analogous to $\frac{1}{\pi}\delta(T)$ in the corresponding breakdown
of $N(T)$ as $\check{N}(T)+S(T)+\frac{1}{\pi}\delta(T)$ for $\zeta$.
The asymptotic piece $\frac{1}{\pi}\delta(T)$ makes a non-trivial
contribution to the root side of the $\mu=-2$ root identity in the
case of $\zeta$.

\textbf{(b)} The form of the divergent piece $\check{N}(T)=\frac{2g}{C}T$
for $\zeta_{k}$ is much simpler than the corresponding divergent
piece $\check{N}(T)=\frac{T}{2\pi}\ln\left(\frac{T}{2\pi}\right)-\frac{T}{2\pi}+\frac{7}{8}$
for $\zeta$, reflecting the much simpler equi-spacing of the generalised
roots for each $\lambda$-factor in the general factorisation of $\zeta_{k}$
in (\ref{eq:zeta_k_general_form2}). In particular, the extra $\ln\left(\frac{T}{2\pi}\right)$
factor in $\check{N}(T)$ for $\zeta$ plays a pivotal role in the
root-side calculations for $r_{\zeta}(s_{0},\mu)$ for $\mu=0,-1$
and $-2$, since it leads to non-trivial differences in Cesaro limiting
behaviour as $T\rightarrow\infty$ and as $\tilde{T}\rightarrow\infty$,
and thereby leads to the non-zero contributions to $r_{\zeta}(s_{0},\mu)$
from the non-trivial roots in these cases. By contrast, for $\zeta_{k}$,
the contributions from the $\check{N}(T)$ pieces as $T\rightarrow\infty$
and as $\tilde{T}\rightarrow\infty$ always cancel fully, leading
to zero overall contributions from this divergent piece to $r_{\zeta_{k}}(s_{0},\mu)$
when $\mu=0,-1$ or $-2$.

(c) For $\zeta_{k}$ the structure of the oscillatory piece, $S(T)$,
is much simpler than it is for the corresponding argument of the Riemann
zeta function in the case of $\zeta$. Consequently, in the case of
$\zeta_{k}$, it is easy to calculate and estimate $S(T)$ and its
integrals, $S_{i}(T)$, in contrast to the situation for $S(T)$ and
the $S_{i}(T)$ for $\zeta$, where we had to rely in {[}1{]} on established
but difficult estimates to obtain Cesaro limits for these functions,
and where these estimates were only unconditional in the cases of
$S$ and $S_{1}$, being conditional on the RH for $\zeta$ for the
higher $S_{i},\: i\geq2$. Note, however, that the contributions from
$S(T)$-related terms to the root sides of the root identities when
$\mu=0,-1$ and $-2$ is zero in both settings (at least modulo the
RH in the case of $\zeta$ when $\mu=-1$ or $-2$).

Overall then, we see from the results of this paper, and from the
comparison between the root identity computations in the setting of
the Riemann zeta function (in {[}1{]}) and of zeta functions of curves
over finite fields here, that it is not merely the existence of both
Euler and Hadamard product formulas that determines the consequences
of these root identities (especially when $\mu=0,-1$ and $-2$).
The detailed, explicit form of the function $N(T)$ is also crucial
in driving the Cesaro computations on the root sides of these identities;
and in turn is thus what drives the apparent contradiction of the
RH in the case of $\zeta$, while remaining consistent with the RH
for $\zeta_{k}$.

\section{Appendices}

\subsection{Tools for performing numerical calculations used in this paper}

The following is the R-code used to generate the example graphs in
section 2, which demonstrate numerically the fact that zeta functions
of curves over finite fields, $\zeta_{k}$, do indeed satisfy the
generalised root identities (result 1). It is straightforward to adapt
this code to conduct further testing as desired.

\paragraph{Code:}

{\small \#\#\#\#\#\#\#\#\#\#\#\#\#\#\#\#\#\#\#\#\#\#\#\#\#\#\#}{\small \par}

\noindent {\small \# set working directory (adapt as appropriate the
argument of the setwd function }{\small \par}

\noindent {\small \# for setting the working directory in the first
line of code)}{\small \par}

\noindent {\small \#\#\#\#\#\#\#\#\#\#\#\#\#\#\#\#\#\#\#\#\#\#\#\#\#\#\#\#}{\small \par}

\noindent {\small setwd(\textquotedbl{}C:/...\textquotedbl{}) }{\small \par}

\noindent {\small options(stringsAsFactors=FALSE)}{\small \par}

\noindent {\small \#\#\#\#\#\#\#\#\#\#\#\#\#\#\#\#\#\#}{\small \par}

\noindent {\small \# Derivative side for given lambda-factor (ignore
nu\_lambda terms as }{\small \par}

\noindent {\small \# same on both sides)}{\small \par}

\noindent {\small \#\#\#\#\#\#\#\#\#\#\#\#\#\#\#\#\#\#}{\small \par}

\noindent {\small DerivSide <- function(s0, mu, q, sigma0, tau0, Nterms)
\{}{\small \par}

\noindent {\small theta0 <- tau0 {*} log(q)}{\small \par}

\noindent {\small lambda <- (q\textasciicircum{}sigma0){*} exp(complex(real=0,imaginary=theta0))}{\small \par}

\noindent {\small p1 <- exp(complex(real=0,imaginary=pi {*} mu))/gamma(mu)}{\small \par}

\noindent {\small f1 <- lambda {*} (log(q))\textasciicircum{}mu {*}
(q\textasciicircum{}(-s0))}{\small \par}

\noindent {\small for (n in 2:Nterms) \{}{\small \par}

\noindent {\small tp <- ((lambda\textasciicircum{}n)/(n\textasciicircum{}(1-mu))){*}((log(q))\textasciicircum{}mu)
{*} (q\textasciicircum{}(-n{*}s0))}{\small \par}

\noindent {\small f1 <- f1+tp}{\small \par}

\noindent {\small \}}{\small \par}

\noindent {\small p1 {*}(sum(f1))}{\small \par}

\noindent {\small \}}{\small \par}

\noindent {\small \#\#\#\#\#\#\#\#\#\#\#\#\#\#\#\#\#\#\#\#\#\#\#\#\#\#\#\#\#\#\#\#\#}{\small \par}

\noindent {\small \# Root side for given lambda factor (ignore nu\_lambda
terms as same }{\small \par}

\noindent {\small \# on both sides)}{\small \par}

\noindent {\small \#\#\#\#\#\#\#\#\#\#\#\#\#\#\#\#\#\#\#\#\#\#\#\#\#\#\#\#\#\#\#\#\#}{\small \par}

\noindent {\small RootSide <- function(s0, mu, q, sigma0, tau0, Nroots)
\{}{\small \par}

\noindent {\small r0\_lambda <- complex(real=sigma0,imaginary=tau0)}{\small \par}

\noindent {\small C <- 2{*}pi/log(q)}{\small \par}

\noindent {\small rootsVec <- 1:Nroots}{\small \par}

\noindent {\small resA1 <- sum(complex(real=rep(s0-sigma0, Nroots),
imaginary=-tau0}{\small \par}

\noindent {\small -C{*}rootsVec)\textasciicircum{}(-mu))}{\small \par}

\noindent {\small resA2 <- sum(complex(real=rep(s0-sigma0, Nroots),
imaginary=-tau0}{\small \par}

\noindent {\small + C{*}rootsVec)\textasciicircum{}(-mu))}{\small \par}

\noindent {\small resA <- resA1+resA2+(complex(real=s0-sigma0, imaginary=-tau0)\textasciicircum{}(-mu))}{\small \par}

\noindent {\small resA<-resA-complex(real=0,imaginary=1){*}(1/((1-mu){*}C)){*}complex(real=s0-sigma0,}{\small \par}

\noindent {\small imaginary=-tau0-C{*}Nroots)\textasciicircum{}(1-mu)
+ complex(real=0,imaginary=1){*}(1/((1-mu){*}C))}{\small \par}

\noindent {\small {*}complex(real=s0-sigma0, imaginary=-tau0+C{*}Nroots)\textasciicircum{}(1-mu)}{\small \par}

\noindent {\small resA <- resA - 0.5{*}complex(real=s0-sigma0, imaginary=-tau0-C{*}Nroots)\textasciicircum{}(-mu) }{\small \par}

\noindent {\small - 0.5{*}complex(real=s0-sigma0, imaginary=-tau0+C{*}Nroots)\textasciicircum{}(-mu)}{\small \par}

\noindent {\small resA<-resA-complex(real=0,imaginary=1){*}mu{*}(C/12){*}complex(real=s0-sigma0,}{\small \par}

\noindent {\small imaginary=-tau0-C{*}Nroots)\textasciicircum{}(-1-mu)
+ complex(real=0,imaginary=1){*}mu{*}(C/12)}{\small \par}

\noindent {\small {*}complex(real=s0-sigma0, imaginary=-tau0+C{*}Nroots)\textasciicircum{}(-1-mu)}{\small \par}

\noindent {\small resA<-resA-complex(real=0,imaginary=1){*}mu{*}(mu+1){*}(mu+2){*}((C\textasciicircum{}3)/720){*}}{\small \par}

\noindent {\small complex(real=s0-sigma0,imaginary=-tau0-C{*}Nroots)\textasciicircum{}(-3-mu)}{\small \par}

\noindent {\small + complex(real=0,imaginary=1){*}mu{*}(mu+1){*}(mu+2){*}((C\textasciicircum{}3)/720)}{\small \par}

\noindent {\small {*}complex(real=s0-sigma0, imaginary=-tau0+C{*}Nroots)\textasciicircum{}(-3-mu)}{\small \par}

\noindent {\small exp(complex(real=0, imaginary=pi{*}mu)){*}resA}{\small \par}

\noindent {\small \}}{\small \par}

\noindent {\small \#\#\#\#\#\#\#\#\#\#\#\#\#\#\#\#\#\#\#\#}{\small \par}

\noindent {\small \# Graph mu > -1.5}{\small \par}

\noindent {\small \#\#\#\#\#\#\#\#\#\#\#\#\#\#\#\#\#\#\#\#}{\small \par}

\noindent {\small s0 <- 5.1238}{\small \par}

\noindent {\small muVec <- seq(-1.45, 2.55, .1)}{\small \par}

\noindent {\small q <- 25}{\small \par}

\noindent {\small sigma0 <- 0.6}{\small \par}

\noindent {\small tau0 <- 3{*}pi/4}{\small \par}

\noindent {\small Nterms <- 20}{\small \par}

\noindent {\small Nroots <- 1000}{\small \par}

\noindent {\small derivVec <- unlist(lapply(muVec, FUN=function(mu)\{DerivSide(s0,
mu, q, sigma0, }{\small \par}

\noindent {\small tau0, Nterms)\}))}{\small \par}

\noindent {\small rootVec <- unlist(lapply(muVec, FUN=function(mu)\{RootSide(s0,
mu, q, sigma0, }{\small \par}

\noindent {\small tau0, Nroots)\}))}{\small \par}

\noindent {\small png(\textquotedbl{}Full Test Re Root and Deriv side
mu geq -1.5 1000.png\textquotedbl{})}{\small \par}

\noindent {\small plot(muVec, Re(derivVec), col=\textquotedbl{}blue\textquotedbl{},
main=\textquotedbl{}Re Root and Derivative side ZetaFF}{\small \par}

\noindent {\small (mu>-1.5)\textquotedbl{}, sub=\textquotedbl{}1000
roots\textquotedbl{})}{\small \par}

\noindent {\small points(muVec, Re(rootVec), col=\textquotedbl{}red\textquotedbl{},
pch=20)}{\small \par}

\noindent {\small legend(\textquotedbl{}topright\textquotedbl{}, fill=c(\textquotedbl{}blue\textquotedbl{},\textquotedbl{}red\textquotedbl{}),
legend=c(\textquotedbl{}Re Deriv side\textquotedbl{}, \textquotedbl{}Re
Root side\textquotedbl{}))}{\small \par}

\noindent {\small abline(h=0)}{\small \par}

\noindent {\small dev.off()}{\small \par}

\noindent {\small png(\textquotedbl{}Full Test Im Root and Deriv side
mu geq -1.5 1000.png\textquotedbl{})}{\small \par}

\noindent {\small plot(muVec, Im(derivVec), col=\textquotedbl{}blue\textquotedbl{},
main=\textquotedbl{}Im Root and Derivative side ZetaFF }{\small \par}

\noindent {\small (mu>-1.5)\textquotedbl{}, sub=\textquotedbl{}1000
roots\textquotedbl{})}{\small \par}

\noindent {\small points(muVec, Im(rootVec), col=\textquotedbl{}red\textquotedbl{},
pch=20)}{\small \par}

\noindent {\small legend(\textquotedbl{}topright\textquotedbl{}, fill=c(\textquotedbl{}blue\textquotedbl{},\textquotedbl{}red\textquotedbl{}),
legend=c(\textquotedbl{}Im Deriv side\textquotedbl{}, \textquotedbl{}Im
Root side\textquotedbl{}))}{\small \par}

\noindent {\small abline(h=0)}{\small \par}

\noindent {\small dev.off()}{\small \par}

\section{Acknowledgements}

We thank Professor Samuel Patterson for suggesting the line of inquiry
pursued in this paper.

\end{document}